\documentclass[11pt,a4paper]{article}

\usepackage{graphicx}
\usepackage[all]{xy}

\newcommand{\arrows}{\,\lower1pt\hbox{$\longrightarrow$}\hskip-.24in\raise2pt
             \hbox{$\longrightarrow$}\,}

\newcommand{\algd}{\mathcal{A}{\rm{lgd}}}

\usepackage{amsmath}
\usepackage{amsthm}
\usepackage{amssymb}

\numberwithin{equation}{section} \allowdisplaybreaks
\newtheorem{theorem}{Theorem}[section]
\newtheorem{lemma}[theorem]{Lemma}
\newtheorem{proposition}[theorem]{Proposition}
\newtheorem{corollary}[theorem]{Corollary}
\newtheorem{definition}[theorem]{Definition}
\newtheorem{example}[theorem]{Example}
\newtheorem{remark}[theorem]{Remark}
\theoremstyle{definition}

\begin{document}
\thispagestyle{empty}

\bigskip

\centerline{\bf Modular classes of Lie algebroid morphisms}

\bigskip

\centerline{Y. Kosmann-Schwarzbach, C. Laurent-Gengoux and A. Weinstein}

\bigskip
\begin{center}
{\textsc{ Dedicated to Bertram
Kostant on his eightieth birthday}}
\end{center}
\bigskip

\noindent{\bf Abstract}
We study the behavior of the modular class of a Lie algebroid
under general Lie algebroid morphisms by introducing the relative 
modular class. 
We investigate the modular classes of pull-back
morphisms and of base-preserving morphisms associated to Lie algebroid
extensions.
We also define generalized morphisms, including Morita
equivalences,
that act on the $1$-cohomology,
and observe that the relative modular class is a coboundary on the
category of Lie algebroids and generalized morphisms with values in
the $1$-cohomology.

\bigskip

\section*{Introduction}
The modular class of a
Poisson manifold first appeared, without a name,
in the work of Koszul \cite{Kz}.
The name was proposed in \cite{W}, where Koszul's class
was shown to be the obstruction to the existence of 
a smooth Poisson trace, and hence the Poisson analogue of 
the modular automorphism group
of a von Neumann algebra.  The extension of the modular class
to Lie algebroids was also proposed there.

This extension was studied in detail by Evens, Lu and Weinstein
\cite{ELW},
who showed that
the modular class of
the cotangent Lie algebroid of a 
Poisson manifold is twice that of the Poisson manifold.
The theory was then developed 
in the framework of Lie-Rinehart algebras by Huebschmann \cite{Hueb},
and the role of Batalin-Vilkovisky 
algebras was emphasized in \cite{Xu} and 
\cite{yks2000}.  

The modular class of a base-preserving morphism of Lie
algebroids was defined in 
\cite{GMM}, using an approach similar to that of \cite{yks2000},
and studied in \cite{KW}, where it is 
called a relative modular class, further explaining results of
\cite{KL} on the modular classes of twisted
Poisson structures on Lie algebroids.

In this paper we extend this relative notion of modular class
to Lie algebroid morphisms which are not necessarily base-preserving.
We look at many examples and study the properties of this new  class.
We introduce a category of generalized morphisms of Lie algebroids
which includes the
weak Morita equivalences of Ginzburg \cite{Gi} along with ordinary
morphisms, and we investigate
the modular class in the context of this category.

In Section  \ref{firstsection}, we recall the definitions of Lie
algebroid morphisms and representations, and we 
 prove that representations
can be pulled back (Proposition \ref{prop1}).
We then recall the
definition of the characteristic class of a Lie algebroid representation 
on a line bundle and that of the
modular class of a Lie algebroid, 
and we study the behavior of the characteristic class of a
representation under  pull-back (Proposition \ref{prop2}).
We can then define the modular class of a Lie
algebroid morphism (Definition \ref{morphism}), an
extension of the definition for
base-preserving morphisms.
In Section \ref{characteristic}, 
the modular class of a Lie algebroid morphism is 
shown to be the characteristic class of a Lie algebroid representation which we
construct explicitly (Theorem \ref{modularmodule}).

Section \ref{thirdsection} is devoted to the modular classes of pull-back Lie
algebroid morphisms; these are the ones for which,
 in some sense, only the base manifold changes. 
We recall the notion of pull-back
 of a Lie algebroid (Definition \ref{pbAlg}) following \cite{M}. 
We first prove that the modular class of a pull-back morphism vanishes when the
morphism covers a submersion (Proposition \ref{submersion}), then we
prove (Theorem \ref{modtransverse})
that this
result is still valid in the more general case of the transverse maps
defined in Definition \ref{transverse}. 
A counter-example (Section \ref{counterex}) shows that the result
is not true in general when the base map
fails to be transverse.
We also show that a surjective submersion induces an isomorphism from the
$1$-cohomology of a Lie algebroid to that of its pull-back (Proposition
\ref{bijective}).

In Section \ref{fourthsection}, we consider some base-preserving Lie 
algebroid morphisms, and we show that the modular class of a Lie algebroid
extension of a transitive
Lie algebroid $B$ by a unimodular
Lie algebroid $C$ can be computed in
terms of a representation of $B$ on a line bundle
(Theorem \ref{th:rel_pb}), generalizing a theorem in \cite{KY}.

In Section \ref{fifthsection}, we define generalized 
Lie algebroid morphisms and the Morita equivalence of Lie algebroids,
as well as the modular class of a generalized morphism with connected and 
simply-connected fibers. We prove that the modular classes of
isomorphic generalized morphisms
are equal, and give a formula for the modular class of the composition of
two generalized morphisms
(Theorem \ref{composition}).
It follows from the definition and from Proposition \ref{submersion} that 
the modular class of a Morita equivalence with connected and simply-connected
fibers vanishes.

In the Appendix (Section \ref{Appendix}), we describe the modular class of morphisms of Lie algebroids 
as a $1$-coboundary on the category of Lie
  algebroids.

There are many interesting questions which remain concerning the
modular classes of Lie algebroids and their morphisms.
\begin{enumerate}

\item
The consequences of the
unimodularity of a Lie algebroid for the existence of a trace on 
an algebra quantizing the Poisson algebra 
of the dual vector bundle have been studied in \cite{NW}.
What is the meaning of the 
modular class 
of a morphism
in terms of these algebras?
  
\item
Modular classes have been shown (\cite{ELW} \cite{Hueb} \cite{Xu})
to be closely related to duality between Lie algebroid homology and
cohomology, and the consequences of the unimodularity of a Poisson
variety for the Van den Bergh dualizing module of its deformation
quantization have been studied in \cite{D}.
What is the meaning of the modular class of a morphism,
in this context?

\item
The modular class is the first in
a sequence of higher
cohomology classes of Lie algebroids which appear in the work of 
Kubarski (see
\cite{Kubarski}), 
Fernandes
\cite{F} and  
Crainic \cite{C}.  How do these higher classes behave under morphisms?

\end{enumerate}

\bigskip

\noindent{\it Conventions.}
All the manifolds, bundles and maps will be
assumed to be smooth.
We denote
the space of (smooth global)
sections of a
real vector bundle $E$ by $\Gamma E$, and its dual vector
bundle by $E ^*$.

A Lie algebroid is denoted by a triple $(A \to M,\rho_A, [~, ~]_A)  $,
where $A \to M$ is a vector bundle, $\rho_A $ is the anchor map, and
$[~, ~]_A $ is the Lie bracket of sections. When no confusion is
possible, we shall abbreviate 
$[~, ~]_A$ to $[~, ~]$.
We often denote a Lie algebroid simply by $A \to M $ or $A$,
in which case the anchor will be understood to be $\rho_A$.

We shall call a topological space $1$-connected if it is both connected and
simply-connected.

\bigskip

\noindent{\it Acknowledgements.}
Alan Weinstein would like to thank the Institut Math\'e\-matique de Jussieu for
hospitality, and the National Science Foundation for partial support
of this research through Grant DMS--0707137.

\section{Morphisms and representations of Lie algebroids}\label{firstsection}
Base-preserving morphisms between Lie algebroids are  simply vector
bundle maps inducing homomorphisms between Lie algebras of sections.
More general morphisms, covering maps between manifolds, were defined
by Higgins and Mackenzie \cite{HM}. Vaintrob
\cite{V} formulated an equivalent condition in terms of vector fields on
supermanifolds. A third formulation, close to Vaintrob's, was given by
Chen and Liu \cite {CL}, who proved the equivalence of their
definition to that of Higgins and Mackenzie. 

\medskip

Modular classes of base-preserving morphisms of Lie
algebroids were defined and studied in \cite{GMM} and \cite{KW}.
Here we review the definition of general
Lie algebroid morphisms, and  we extend the notion of modular
class to them.

\subsection{Morphisms of Lie algebroids}
Let $A \to M$ and $B \to N$ be
real vector bundles, and let $ (\Phi, \phi)$ be a vector bundle map,
 $$ \begin{array}{ccc}  A &\stackrel{\Phi}{\to} &B  \\ \downarrow &
 &\downarrow  \\ M & \stackrel{\phi}{\to} & N \\
 \end{array}$$
Such a map induces a pull-back operator $\widetilde \Phi^*:\Gamma (B^*) \to \Gamma (A^*)$ on sections, 
defined by
\begin{equation}
< (\widetilde \Phi^* \beta)_m, a_m > \, = \,  < (\beta \circ \phi)(m),
\Phi(a_m) > \ ,
\end{equation}
for $\beta \in \Gamma(B^*)$, $m \in M$ and $a \in
\Gamma A$.
This operator is induced by the
base-preserving morphism of vector bundles
over $M$, $\Phi^* \colon \phi^{!}B^* \to A^*$, where
$\phi^{!}B^*$ is the pull-back of $B^*$ under $\phi$.
We may extend $\widetilde \Phi^*$ to an exterior algebra
homomorphism, $\wedge ^{\bullet} \widetilde \Phi^*:
\Gamma (\wedge^{\bullet} B^*) \to \Gamma (\wedge^{\bullet} A^*)$,
where
we set 
$(\wedge ^{0} \widetilde \Phi^*)(f)$ 
to be $f \circ \phi$,
for $f \in C^\infty(M)$.
Recall that,
for any Lie algebroid $(A \to M,\rho_A,[~,~]_A) $,
the Lie algebroid differential $d_A$ is defined by
$$
\begin{array}{cc}
(d_A \alpha)(a_0, a_1, \ldots, a_k)  = \sum_{i=0}^{k} 
(-1)^i
\rho_A(a_i) \cdot \alpha(a_0, \ldots, \widehat{a_i}, \ldots, a_k) \\
 +  \sum_{0 \leq i < j \leq k} 
(-1)^{i+j}
 \alpha([a_i,a_j]_A, a_0, \ldots, \widehat{a_i}, \ldots,
 \widehat{a_j}, \ldots, a_k) \ ,
\end{array}
$$
 for $\alpha \in
\Gamma (\wedge^{k}A^*)$, $k \in \mathbb N$, $a_0, \ldots, a_k \in
\Gamma A$. The differential $d_A$ turns $\Gamma (\wedge^{\bullet}
A^*)$ into a complex, whose cohomology is called the \emph{Lie
  algebroid cohomology}
and is denoted by $H^\bullet(A)$.  
It is natural to make the following definition.
\begin{definition}\label{def1}
Let $A \to M$  and $B \to N$
be Lie algebroids, and let $(\Phi,\phi)$ be a vector bundle map from
$A$ to $B$. The pair $(\Phi,\phi)$ is a {\em morphism of Lie
  algebroids}
if the map $\wedge^{\bullet} \widetilde \Phi^* :
(\Gamma(\wedge^{\bullet}B^*), d_B) \to
(\Gamma(\wedge^{\bullet}A^*), d_A)$ is a chain map.
\end{definition}

When $(\Phi,\phi)$ is base-preserving, {\it{i.e.}}, $\phi$ is an identity map,
 $\widetilde \Phi^*$ is the usual dual of
$\Phi$, and the condition that
$\wedge^{\bullet} \widetilde \Phi^*$ be a chain map is equivalent to
the usual definition 
of a morphism of Lie algebroids as a vector bundle map preserving
anchors and brackets \cite{M1} \cite{M}.  
In fact, Definition \ref{def1} is just the ``ordinary language'' version of the 
definition due to Vaintrob \cite{V}:
$(\Phi, \phi)$ is a morphism if the homological vector
fields $d_A$ and $d_B$ are $\Pi\Phi$-related, where $\Pi \Phi$ is
$\Phi$ considered as a morphism  from $\Pi A$ to 
$\Pi B$, the supermanifolds obtained from $A$ and $B$ by ``making the
fibers odd,'' {\it {i.e.}}, where the
functions on $\Pi A$ are the sections of $\wedge ^{\bullet} A^*$,
and the functions
on $\Pi B$ are the sections of $\wedge ^{\bullet} B^*$.

To summarize this section, we may say that the morphisms make Lie
algebroids over all manifolds a category
$\mathcal{A}{\rm{lgd}}$. 
The Lie algebroid complex is a contravariant
functor from this category to the category of complexes of vector
spaces over $\mathbb R$, and Lie algebroid cohomology is a contravariant
functor from 
$\mathcal{A}{\rm{lgd}}$ to the category of graded
real vector spaces.  We discuss this categorical approach further in
the Appendix.

\subsection{Representations of Lie algebroids}\label{basic}

It is useful to think of $\mathcal{A}{\rm{lgd}}$ as an enlargement of
the category of manifolds and smooth maps, which form a full
subcategory when we identify each manifold $M$ with its tangent bundle
Lie algebroid.  The restriction of
Lie algebroid cohomology is, of course, just de Rham cohomology.
Other full subcategories are that of 
Lie algebras, for which
the Lie algebroid cohomology is the
Chevalley-Eilenberg cohomology, and that of the zero Lie algebroids.

The theory of flat connections on vector bundles over manifolds
becomes a representation theory when we pass from manifolds to Lie
algebroids, 
and this representation theory generalizes that of Lie algebras.

We recall that a {\it representation} \cite{M1} \cite{M} 
of a Lie algebroid $A$
on a vector bundle $E$, both over the base $M$,
is an $\mathbb R$-bilinear map
$ D : \Gamma A \times \Gamma  E  \to \Gamma E $,
denoted $ (a,x) \mapsto D(a) x $, or simply $ a \cdot x$ when no
confusion is possible,
such that, for all $f \in
C^{\infty}(M)$, $a, a_1, a_2 \in \Gamma A$ and $x \in \Gamma E$,
 $$ \left\{ \begin{array}{lrcl} (i) & (fa) \cdot x  &= & f (a \cdot x)
     \ ,\\
(ii) & a  \cdot (fx) &= &  f (a \cdot x) +  (\rho_A a)  (f) \,  x \ , \\
(iii) &  a_1 \cdot (a_2 \cdot x) -  a_2 \cdot (a_1 \cdot x)  &=&
[a_1,a_2]_A
\cdot x \ . \\
 \end{array} \right. $$
Equivalently, a {\it representation of $A$ on $E$}, also called a {\it flat
$A$-connection on $E$} or an {\it $A$-module structure on $E$}, is a Lie
algebroid morphism over $M$ 
from $A$ to $\mathcal D E$, the Lie algebroid of
derivations on $E$ whose sections are the covariant differential
operators (CDO's) or derivative
endomorphisms of $\Gamma E$ (see \cite{M}). 
More generally, an
{\it $A$-connection on $E$} is a vector bundle morphism over $M$ 
from $A$ to $\mathcal D E$, {\it {i.e.}}, a map satisfying (i) and (ii)
above. We also recall the following constructions.
\begin{itemize}
\item If $D$ is a representation of $A$ on $E$, then the dual
  representation $D^*$
is the representation of $A$ on $E^*$ defined by
$$
< D^*(a) \xi , x > \, = \, - < \xi, D(a) x > + \, (\rho_A a)(<\xi,x>) \ ,
$$
for all $a \in \Gamma A$, $x \in \Gamma E$, $\xi \in \Gamma (E^*)$.
\item If $D_1
$ and $D_2$ are representations of $A$ on 
vector bundles over $M$,
$E_1$ and
  $E_2$, then $D_1 \otimes D_2$ is the representation of $A$ on
  $E_1 \otimes E_2$ defined by
$$
(D_1 \otimes D_2)(a) (e_1 \otimes e_2) = D_1(a) (e_1) \otimes e_2 +
e_1 \otimes D_2 (a) (e_2) \ ,
$$
for all $a \in \Gamma A$,
 $e_1 \in \Gamma (E_1)$ and $e_2 \in \Gamma (E_2)$.
\end{itemize}

We call an
$\mathbb R$-linear endomorphism of $\Gamma (\wedge^{\bullet} A^*
\otimes E)$ of degree~$1$
and of square zero
a {\it differential} on $\Gamma (\wedge^{\bullet} A^*
\otimes E)$.
The following generalization of a well-known characterization of
flat connections  
proves that the definition of a representation of a Lie
algebroid $A$ adopted above is equivalent to that of an
$A$-module in \cite{V}.

\begin{proposition}\label{differential}
Any
differential $d_{A,E}$ on $\Gamma (\wedge^{\bullet} A^*
\otimes E)$
satisfying \begin{equation}\label{leibniz2}
d_{A,E}(\alpha \otimes x) = d_A \alpha \otimes x + (-1)^{|\alpha|} \alpha
\otimes d_{A,E} \, x \ ,
\end{equation}
for all $\alpha \in \Gamma (\wedge^{\bullet} A^*)$ of degree~$|\alpha|$, and $x
\in \Gamma E$,
gives rise to a representation of $A$ on $E$
defined by
\begin{equation}
a \cdot x = \iota_a(d_{A,E} \, x) \ ,
\end{equation}
for $a \in  \Gamma A$ and $x \in \Gamma E$, where $\iota$ denotes 
the interior product.
Conversely, each representation of $A$ on $E$ gives rise in this way to
a differential on $\Gamma (\wedge^{\bullet} A^*
\otimes E)$ satisfying (\ref{leibniz2}).
\end{proposition}
\begin{proof}
The
$\mathbb R$-linear map
${d_{A,E}}_{| \Gamma E} \colon \Gamma E \to \Gamma(A^* \otimes E)$ satisfies,
for each\break $x \in \Gamma E$ and $f \in C^{\infty}(M)$,
\begin{equation}\label{leibniz}
d_{A,E}(fx) = f d_{A,E} x +  d_{A}f \otimes x \ ,
\end{equation}
and therefore $(a,x) \mapsto a \cdot x$ satisfies (i) and (ii)
above.

Conversely, given a representation, the map $d_{A,E}$ is
well-defined on
 $\Gamma E$
because of (i) and satisfies \eqref{leibniz} because
of (ii). It can then be uniquely extended to an
$\mathbb R$-linear endomorphism of degree~$1$ 
of
$\Gamma(\wedge^{\bullet}A^* \otimes E)$ satisfying \eqref{leibniz2}.

We claim that $(d_{A,E})^2 = 0$ is equivalent to the flatness property (iii).
In fact, expressing 
$d_{A,E} \, x$, for $x \in \Gamma E$,
locally as a finite
sum $d_{A,E} \, x = \sum_{k} \alpha_k \otimes y_k$,
with $\alpha_k \in \Gamma (A^*)$
and $y_k \in \Gamma E$, yields, for $a$ and $a' \in \Gamma A$,
$$
\begin{array}{lll}
\iota_{a'} \iota_{a} (d_{A,E}(d_{A,E} \, x))  =  \sum_k
\iota_{a'} \iota_a (d_A \alpha_k \otimes y_k - \alpha_k \otimes d_{A,E} \,
y_k) \\
=  \sum_k \left ((d_A \alpha_k)(a,a') y_k + < \alpha_k , a' > \iota_a d_{A,E} y_k -
 < \alpha_k , a > \iota_{a'} d_{A,E} \, y_k \right ) \ ,
\end{array}
$$
while
$$
\begin{array}{lll}
a \cdot (a' \cdot x) - a' \cdot (a \cdot x) 
&=& \sum_k  ( \iota_{a}
d_A <\alpha_k,a'> y_k
- \iota_{a'} d_A <\alpha_k,a> y_k\\
&+& <\alpha_k,a'> \iota_{a} d_{A,E} \, y_k
-<\alpha_k,a> \iota_{a'} d_{A,E} \, y_k ) \ .
\end{array}
$$
Therefore
$$
a \cdot (a' \cdot x) - a' \cdot (a \cdot x) - [a,a'] \cdot x =
\iota_{a'} \iota_a (d_{A,E} (d_{A,E} \, x)) \ ,
$$
which proves the claim.
\end{proof}

This proof shows that, more generally, 
the $A$-connections on $E$ are in one-to-one
correspondence with the  $\mathbb R$-linear endomorphisms of degree
$1$ of $\Gamma (\wedge^{\bullet} A^* \otimes E)$
satisfying (\ref{leibniz2}).

If $A \to M$ and $B \to N$ are Lie algebroids 
and $D_A$ (resp., $D_B$) is a representation of 
$A \to M$ (resp., $B \to N$) on a
vector bundle $E \to M$ (resp., $F \to N$), a vector bundle map $\Psi:E \to F$ is
said to be a {\em morphism of representations from $D_A$ to $D_B$ covering the 
Lie algebroid morphism} $(\Phi, \phi)$ 
from $A \to M$ to $B \to N$ if $\Psi: E \to F $ covers $\phi: M \to N $,
and if the map 
$\wedge^{\bullet} \widetilde \Phi ^* \otimes \widetilde \Psi^*$ is a
chain map from the complex $(\Gamma(\wedge^{\bullet} B^* \otimes F^*),
d_{B,F^*})$ to the complex
$(\Gamma(\wedge^{\bullet} A^* \otimes E^*),d_{A,E^*})$, {\it{i.e.}}, if the
following diagram commutes.
\begin{equation}\label{chain0} \xymatrix{
  \Gamma(\wedge^{\bullet} B^* \otimes F^*) \ar[r]^{d_{B,F^*}}
  \ar@<0.5ex>[d]^{\wedge^{\bullet}{\widetilde \Phi^*} \otimes
    \widetilde \Psi^* }  &
  \Gamma(\wedge^{\bullet + 1}
 B^* \otimes F^*) \ar@<0.5ex>[d]^{\wedge^{\bullet + 1}{\widetilde \Phi^*}
   \otimes \widetilde \Psi^* } \\
  \Gamma(\wedge^{\bullet} A^* \otimes E^*)  \ar[r]^{d_{A,E^*}} &
  \Gamma(\wedge^{\bullet + 1} A^* \otimes E^*) }
\end{equation}

If $M=N=\{{\rm{pt}}\}$, the commutativity of \eqref{chain0} reduces to the
condition that $\Psi$ be a morphism from the $A$-module $E$ to the
$B$-module $F$ covering the Lie algebra morphism $\Phi:A \to B$,
$\Psi(D_A(a \otimes v)) = D_B(\Phi(a) \otimes \Psi(v))$, for all $a
\in A$ and $v \in E$.

\subsection{Pull-back of a representation by a Lie algebroid
  morphism}\label{pullbackrep}
Let $F \to N$ be a vector bundle, and let
$\phi^! F \to M $ be its pull-back by a map
$\phi : M\to N$.  
Any section $y \in \Gamma F $
pulls back to 
the
section $ \phi^!{y} \in \Gamma(\phi^! F)$, where
$\phi^! y = y \circ \phi$. 
The next
proposition generalizes the facts that flat bundles pull back to flat
bundles and that  Lie algebra
representations pull back to representations.

\begin{proposition}\label{prop1}
If $D$ is a representation of the Lie algebroid $B\to N$ on
 $F \to N$, and if $(\Phi, \phi) $ is
a Lie algebroid morphism from $ A \to M $ to $B
\to N $, then there exists a unique representation $\Phi^! D$ of  $A$ on
$\phi^! F$ such that, for all $a \in \Gamma A$, $b \in  \Gamma B$
satisfying
$\Phi \circ a = b \circ \phi$, and $y \in \Gamma F$,
\begin{equation}\label{pullback}
a \cdot (\phi^! y) = \phi^!(b \cdot y) \ .
\end{equation}
Explicitly,
$(\Phi^! D)(a)(\phi^!y) = \phi^! (D(b)y)$.
\end{proposition}
\begin{proof}
We shall define the representation $\Phi^! D$ by means of the
associated differential $d_{A, \phi ^! F} \colon
\Gamma(\phi ^! F) \to \Gamma(A^* \otimes  \phi ^! F)$
described in Proposition \ref{differential}.
We first define the map $d_{A, \phi ^! F}$
on the sections of $\phi ^! F$ of the form $\phi^! y$, for $y \in
\Gamma F$, by
\begin{equation}
d_{A, \phi ^! F}(\phi ^! y) = ( \widetilde \Phi^* \otimes
\phi^!)(d_{B,F} \, y) \ .
\end{equation}
For $f\in C^{\infty}(M)$, we set
\begin{equation}\label{leibniz3}
d_{A, \phi ^! F} (f \otimes \phi^! y)= f \otimes d_{A, \phi ^! F} (\phi^! y)
+  d_{A} f  \otimes \phi^! y \ ,
\end{equation}
A computation that uses the fact that
$(\Phi, \phi)$ is a morphism of Lie algebroids, which implies that
$ \widetilde \Phi^*(d_B \, h) = d_A (h\circ \phi)$ for $h \in C^{\infty}(N)$,
shows that
$$
d_{A, \phi ^! F}((h \circ \phi) f \otimes_{\mathbb R} \phi^!
y) =
d_{A, \phi ^! F}(f \otimes_{\mathbb R} \phi^! (hy)) \ .
$$
Therefore, there is a well-defined
linear map $d_{A, \phi ^! F}$ 
to $\Gamma(A^* \otimes \phi^!F)$ from
$(C^{\infty}(M) \otimes_{\mathbb R} \Gamma F)/ \phi^!(C^{\infty}(N))$,
which is isomorphic to $\Gamma(\phi^!F)$ as a $C^{\infty}(M)$-module.
Because of \eqref{leibniz3}, this map can be extended to
an $\mathbb
R$-linear endomorphism
of degree 1 of $\Gamma( \wedge^{\bullet} A^* \otimes
\phi^!F)$ satisfying \eqref{leibniz2}.
A straightforward computation shows that this
map is of square
zero.
The corresponding representation satisfies, for $y \in \Gamma F$,
$$
a \cdot (\phi^! y)
= \iota_a d_{A, \phi ^! F} (\phi ^! y)= \iota_a
(\widetilde \Phi^* \otimes
\phi^!)(d_{B,F} \, y) \ .
$$
Assuming, without loss of generality, that
$d_{B,F} \, y = \beta \otimes z$, with $\beta \in \Gamma B^*$ and $z \in
\Gamma F$, and assuming that $\Phi \circ a = b \circ \phi$, we obtain
$$
a \cdot (\phi^! y) = <\widetilde \Phi ^* \beta, a > \phi^ ! z
= (< \beta, b > \circ \, \phi) \phi ^! z = \phi ^! (<\beta,b> z)
= \phi ^! (b \cdot y) \ ,
$$
as required. This property defines the pull-back representation uniquely.
\end{proof}

See Remark \ref{remarksubmersion} below for a simple proof of this
proposition in the case where
$\phi$ is a transverse map.

When $d_{B,F}$ and $d_{A,\phi^! F}$ are 
the differentials associated by Proposition \ref{differential} to 
a representation $D$ of $ B \to N $ on $F$ and the pull-back
representation $\Phi^!D$
of $A \to M$ on $\phi^! F$,
the following diagram commutes.
\begin{equation}\label{chain1}
 \xymatrix{
\Gamma(\wedge^{\bullet} B^* \otimes F) \ar[r]^{d_{B,F}}
\ar@<0.5ex>[d]^{\wedge^{\bullet}{\widetilde \Phi}^* \otimes \phi^{!} }
& \Gamma(\wedge^{\bullet + 1 }
B^* \otimes F) \ar@<0.5ex>[d]^{\wedge^{\bullet + 1}{\widetilde \Phi}^* \otimes \phi^! } \\
\Gamma(\wedge^{\bullet} A^* \otimes \phi^{!} F)
\ar[r]^{d_{A,\phi^{!}F}} 
& \Gamma(\wedge^{\bullet + 1} A^* \otimes \phi^{!}F) } 
\end{equation}

\begin{corollary}\label{chainmap}
The canonical projection $\phi_F \colon \phi^! F \to F$
 is a morphism of representations from $D$ to $\Phi^!D$
covering the Lie algebroid morphism $(\Phi, \phi)$.
\end{corollary}
\begin{proof}
Let $ \Phi^{!} (D^*) $ be the pull-back by $(\Phi,\phi) $
of the representation $D^*$ 
of $B \to N $ on $F^*$, which is a 
representation of $A \to M$ on $\phi^! (F^*)$.
By (\ref{chain1}), the following diagram 
commutes.
$$ \xymatrix{
  \Gamma(\wedge^{\bullet} B^* \otimes F^*) \ar[r]^{d_{B,F^*}}
  \ar@<0.5ex>[d]^{\wedge^{\bullet}{\widetilde \Phi^*} \otimes \phi^{!} }  &
  \Gamma(\wedge^{\bullet + 1}
 B^* \otimes F^*) \ar@<0.5ex>[d]^{\wedge^{\bullet + 1 }{\widetilde
     \Phi^*} \otimes \phi^! } \\
  \Gamma(\wedge^{\bullet} 
A^* \otimes \phi^{!}(F^*))  
\ar[r]^{d_{A,\phi^{!}(F^*)}} & \, \, \, \,  \Gamma(\wedge^{\bullet + 1} A^* \otimes
\phi^{!}(F^*)) }$$
When 
$\phi^! (F^*) $ is identified 
to $ (\phi^! F)^*$ as vector bundles, the
map $\phi^! : \Gamma({F^*}) \to \Gamma(\phi^!({F^*}))$ 
is identified with $(\widetilde {\phi_F})^* : \Gamma({F^*}) \to
\Gamma((\phi^! {F})^*)$,
and the previous commutative diagram becomes:
$$ \xymatrix{
  \Gamma (\wedge^{\bullet} B^* \otimes F^*) \ar[r]^{d_{B,F^*}}
  \ar@<0.5ex>[d]^{\wedge^{\bullet}{\widetilde \Phi^*} 
\otimes (\widetilde {\phi_F})^* }  &
  \Gamma (\wedge^{\bullet + 1} B^* \otimes F^*)
  \ar@<0.5ex>[d]^{\wedge^{\bullet + 1}{\widetilde \Phi^*} \otimes
    (\widetilde {\phi_F})^* } \\
  \Gamma (\wedge^{\bullet} A^* \otimes (\phi^{!}F)^*)
  \ar[r]^{d_{A,(\phi^{!}F)^*}} & \, \, \Gamma(\wedge^{\bullet + 1} A^* \otimes
  (\phi^{!} F)^*) }$$
By (\ref{chain0}), the claim is proved. 
\end{proof}

\section{Modular classes of morphisms}\label{sectionmodularclasses}

\subsection{Modular classes of Lie algebroids}
\label{char}
We recall some results of \cite{ELW}.
When $D$ is a representation of the Lie algebroid $A \to M$ on an orientable
line bundle
$L \to M$, and when  $\lambda$ is a nowhere-vanishing section of $L$,
the section $\alpha_{\lambda}$ of $A^*$
defined by
$$
<\alpha_{\lambda}, a > \, \lambda = a \cdot \lambda \ ,
$$
for all $a \in \Gamma A$,
is $d_A$-closed. As above, we have denoted $D(a)
\lambda$ by $a \cdot \lambda$.\break 
We call $\alpha_\lambda$
the {\it characteristic cocycle} associated to the representation $D$
and the section $\lambda$.  
Its class in the $1$-cohomology of $A$ is
independent of the choice of $\lambda$.   We call it
 the {\it characteristic class} of the representation
$D$ and denote it by ${\rm{char}} \, D$.
When $L$ is not orientable, the characteristic class is defined as
one-half that of the representation $D \otimes D$ in the square of the
line bundle $L$.  This is consistent with part (ii) of the following
proposition, which establishes some basic properties of 
characteristic classes. 

\begin{proposition}\label{prop2}
(i)
If $D^*$ is the dual of $D$, then ${\rm{char}} \, (D^*) = -
  \,
{\rm{char}} \, D$.

\noindent(ii)
If $D_1$ and $D_2$ are representations of $A$ on line bundles
$L_1$ and $L_2$ over $M$, then $${\rm{char}} \, ({D_1
  \otimes D_2}) = {\rm{char}} \, ({D_1}) + {\rm{char}} \, ({D_2}) \ . $$

\noindent(iii)
If $(\Phi, \phi)$ 
is a morphism of Lie algebroids from $A \to M$
to $B \to N$, and $D$ is a representation of $B$ on a line bundle,
then
\begin{equation}\label{pullbackchar}
{\rm{char}} \, ({\Phi^!D}) =  \widetilde \Phi^*
({\rm{char}} \, {D}) \ .
\end{equation}
\end{proposition}
\begin{proof} The proofs of (i) and (ii) are straightforward. Let us
  prove (iii).
By Proposition \ref{prop1}, if $D$ is a representation of $B$
on the line bundle $L$, $\Phi^!D$ is a representation of $A$ on
$\phi^{!}L$.
We denote by $\beta$ the representative of ${\rm{char}} \, D$ associated to
$\nu$, a nowhere-vanishing section of $L$, and  by $\alpha$ the
representative of ${\rm{char}} \, ({\Phi^! D})$ associated to $\phi^!
\nu$.
By definition, for all $a \in \Gamma A$,
$$
< \alpha, a> \phi^! \nu = a \cdot (\phi^! \nu)
= \iota_a(d_{A, \phi^! L}(\phi^! \nu)) \ .
$$
Since diagram (\ref{chain1}) commutes,
$$
d_{A, \phi^! L}(\phi^! \nu) = (\widetilde \Phi ^* \otimes \phi
^!)(d_{B,L} \nu) \ .
$$
Therefore, for $a \in \Gamma A$, $\alpha \in \Gamma(A^*)$,
at any $m \in M$,
$$
< \alpha_m, a_m> (\phi^! \nu)_m =
\iota_{a_m}(( \widetilde \Phi ^* \otimes \phi
^!)(d_{B,L} \nu)) (m)
= \iota_{\Phi(a_m)}(d_{B,L}
\nu)(\phi(m))
$$
$$
= <\beta_{\phi(m)}, \Phi(a_m)> \nu_{\phi(m)} = <(\widetilde \Phi ^* \beta)_m, a_m> (\phi^!\nu)_m \ .
$$
We have proved that $\alpha = \widetilde
{\Phi}^* \beta$, and \eqref{pullbackchar} follows.
\end{proof}

Each Lie algebroid $A\to M$
has a canonical representation $D^A$ in the line
bundle 
\begin{equation}\label{linebundle}
L^A = \wedge^{\rm {top}} A \otimes \wedge^{\rm {top}} (T^*M)
\ ,
\end{equation}
defined by
\begin{equation}\label{repA}
D^A_a(\omega \otimes \mu) = [a, \omega]_A \otimes \mu + \omega
\otimes {\mathcal L}_{\rho_A a} \, \mu \ ,
\end{equation}
for all $a \in \Gamma A$, where $\omega \otimes \mu$ is a
nowhere-vanishing section of $L^A$. Here $[~,~]_A$ is the
Gerstenhaber bracket on $\Gamma(\wedge^{\bullet}A)$, while
$\mathcal L$ denotes the Lie derivation.
A section $\alpha$ of $A^*$ satisfying 
\begin{equation}
<\alpha, a> \omega \otimes \mu = D^A_a (\omega \otimes \mu) \ ,
\end{equation}
for all $a \in \Gamma A$, is called the {\it modular cocycle for} $A$
associated to the 
nowhere-vanishing section $\omega \otimes \mu$  
of $L^A$, and the characteristic class it defines 
is called the {\it modular class of} $A$.

\subsection{Definition of the modular class of a morphism}
Using the cohomology pull-back operation associated to any
morphism of Lie algebroids, we may make the following definition.
\begin{definition}\label{morphism}
Let ${\rm{Mod}} \, A$ and ${\rm{Mod}} \, B$ be the modular classes of Lie
algebroids $A \to M $ and $B \to N$, and let $(\Phi, \phi)$  be
a morphism of Lie algebroids from $A \to M$ to $B \to N$.
The {\em relative modular class} or simply the {\em modular class}
of $(\Phi, \phi)$  is
the class, $
{\rm{Mod}} \, {\Phi}$, in the {\rm 1}-cohomology of $A$ defined by
$$
{\rm{Mod}} \, {\Phi} = {\rm{Mod}} \, A -  \widetilde \Phi^*
({\rm{Mod}} \, B) \ .
$$
\end{definition}
This definition extends to general morphims that of the modular class
for base-preserving morphisms given in \cite{GMM} \cite{KW}.

It follows from the definition that, for Lie algebroids $A \to M$, $B
\to N$ and $C \to R$, and morphisms $(\Phi, \phi)$ from $A \to M$ to $B \to
N$ and $(\Psi, \psi)$ from $B \to N$ to $C \to R$,
\begin{equation}\label{comp}
{\rm{Mod}} (\Psi \circ {\Phi}) = {\rm{Mod}} \, \Phi +  \widetilde \Phi^*
({\rm{Mod}} \, \Psi) \ .
\end{equation}
In the particular case of a base-preserving morphism,
${\widetilde{\Phi}}^*$ reduces to $\Phi^*$. Therefore 
(\ref{comp})
(see also (\ref{eq3})) 
is the generalization of relation (3)
of \cite{KW} to the
Lie algebroid morphisms which are
not necessarily base-preserving.

In the following sections, we shall prove that these relative modular classes
are characteristic classes of representations,
and we shall study various examples of such classes.

\subsection{Modular classes of morphisms as characteristic classes}
\label{characteristic}
We will now show that the modular class of a Lie algebroid morphism
from $A \to M$ to $B \to N$ is the characteristic class
of a representation of $A$ in a line bundle, thus generalizing Theorem
3.3 of \cite{KW} to the case of morphisms which are not necessarily
base-preserving.

For any Lie algebroid morphism  $(\Phi, \phi)$ 
from $A \to M$ to $B \to N$, we set
$$L^{\phi} = L^A \otimes \phi ^! ((L^B)^*) \ ,
$$
and
$$
D^{\Phi} = D^A \otimes 
\Phi^!
((D^B)^*) \ ,
$$
where $L^A$ and $D^A$ are defined by (\ref{linebundle}) and
(\ref{repA}), and $\Phi ^! ((D^B)^*)$ is defined in Proposition \ref{prop1}.
\begin{theorem}\label{modularmodule}
Let $(\Phi, \phi)$ be a Lie algebroid morphism  from $A \to
M$ to $B \to N$. 

\noindent(i) $D^{\Phi}$ is a representation of $A$ on
$L^{\phi}$.

\noindent (ii) The modular class, ${\rm {Mod}} \, {\Phi}$, of the morphism
$(\Phi, \phi)$ is the characteristic class of the representation $D^{\Phi}$.
\end{theorem}
\begin{proof}
The fact that $D^{\Phi}$ is  a representation follows from the basic
properties stated in Sections \ref{basic} and
\ref{pullbackrep}, while (ii) follows from 
parts (i) and (iii) of 
Proposition \ref{prop2}.
\end{proof}

The following simple examples will be useful in the next sections.

\begin{example}\label{ex:restriction}\rm{
Let $A \to M $ be a Lie algebroid, $U \subset M $ an open subset
and $A_{|_U} \to U$  the restriction of $A \to M $ to $U$.
Then the modular class of the inclusion 
morphism of $A_{|_U} \to U $ into $A \to M $ vanishes.}
\end{example}

\begin{example}\label{ex:isomorphism}
\rm{The modular class of a Lie algebroid isomorphism vanishes.
In fact, a Lie algebroid 
isomorphism $(\Phi, \phi) $ from $A \to M $ to $B \to N $
induces an isomorphism from the representation of $ A \to M$ on the line
bundle 
$L^A $ to the representation of 
$B \to N$
on the line bundle $L^{B}$
covering $(\Phi,\phi)$.} 
\end{example}

\section{Pull-back Lie algebroids and pull-back morphisms}\label{thirdsection}
\subsection{Definition of pull-backs}

As in Section \ref{pullbackrep}, we denote 
by $\phi^! B \to M $ the pull-back
of a vector bundle $B \to N $
by a map $\phi: M \to N $,
and for any section $b \in \Gamma B $,
we denote by $ \phi^!{b} $ the 
pulled-back 
section of $\phi^! B \to M . $

\begin{definition}\label{admissible}
Let $B \to N$ be a Lie algebroid with anchor $\rho_B$.
A map $\phi \colon M \to N$ is called {\emph{admissible}} if the pull-back
$B \oplus_{TN} TM$, whose fiber at $m \in M$ is $\{ (b,u) \in 
B_{\phi (m)} \oplus T_m M \, | \, 
\rho_B b = (T\phi) u \}$, has a rank independent of $m$, in which case
it is a vector sub-bundle of $\phi^! B \oplus TM$.
\end{definition}
For instance, any surjective submersion is admissible. If $B \to N$ is 
transitive, {\it{i.e.}}, if
$\rho_B$ is
surjective, any map into $N$ is admissible. The inclusion $\mathcal O
\to N$, an
injective immersion, of any orbit $\mathcal O$ of $B\to N$,
{\it{i.e.}}, integral manifold of the distribution on $M$ defined by
the image of $\rho_B$, is admissible.

The following definition is due to Higgins and Mackenzie \cite{HM}
\cite{M}.

\begin{definition}\label{pbAlg}
The pull-back of the Lie algebroid
$(B \to N, \rho_B, [~,~]_B)$
by an admissible map $\phi: M \to N$ is the Lie algebroid
$(\phi^{!!} B \to M, \rho^{!!}_B, [~,~]_{\phi^{!!} B})$,
where
\begin{enumerate}
\item the total space of $\phi^{!!}B \to M$ is $B \oplus_{TN} TM$,
\item the anchor map $\rho^{!!}_B$ is the projection onto the second
  component, and
\item the bracket is defined as follows. For $u, v \in \Gamma (TM)$,
$b_i, c_i \in \Gamma B$, $f_i, g_i \in C^\infty(M)$,
satisfying $\sum_i (f_i \circ \phi)(m)\rho_B(b_i(\phi(m)) =
(T_m\phi) (u(m))$ and 
$\sum_i (g_i \circ \phi)(m)\rho_B(c_i(\phi(m)) =
(T_m\phi) (v(m))$, for all $m \in M$,
$$
\begin{array}{lll}
 [({\sum_i f_i \otimes b_i},u), ( {\sum_i g_i \otimes c_i}, v)]_{\phi^{!!}B}\\
 = \left( \sum_{i.j} f_i g_j \otimes [b_i,c_j]_B + \sum_i (u \cdot g_i)
\otimes c_i  - \sum_i (v \cdot f_i)
\otimes b_i , [u,v]_{TM} \right) \ .
\end{array}
$$
\end{enumerate}
\end{definition}

For example, the pull-back of $B\to N$ by the inclusion $\mathcal O \to
N$ of an orbit is a transitive Lie algebroid over
$\mathcal O$ (see \cite{M}).  

 \begin{remark}
{\rm
In the case when $\phi$ is a submersion,
the bracket, $[~,~]_{\phi^{!!}B}$, is the unique Lie algebroid bracket
such that for any
$b, c  \in \Gamma B $ and $u, v 
\in \Gamma(TM)$
satisfying 
$\rho_B(b) = (T\phi) (u)$ and $\rho_B(c) = (T\phi) (v)$, 
\begin{equation} \label{bracketforsubmersions}
 [(\phi^!{b},u) ,(\phi^!{c},v) ]_{\phi^{!!}B}=
  (\phi^!{[b,c]_B},[u,v]_{TM})  \ .
\end{equation}
}
\end{remark}

The projection onto the first component
is a morphism of Lie algebroids from  $\phi^{!!}B \to M $
to $B \to N $, which we denote by $\phi^{!!}_B $.

\medskip

As an example, we shall determine the pull-back Lie algebroid of
the Lie algebroid $\mathcal D F$ of derivations on a vector bundle
$F$, which exists since $\mathcal D F$ is transitive.
\begin{proposition}\label{derivations}
For any vector bundle $F$ over $N$,
and map $\phi \colon M \to N$,
the Lie algebroids $\phi^{!!}({\mathcal D} F)$ and
$\mathcal D (\phi^! F)$ are isomorphic.
\end{proposition}
\begin{proof}
Recall that a representation  of a Lie algebroid $A \to M $ on a
vector bundle $E \to M $
can be seen as a Lie algebroid morphism from $A \to M$ to ${\mathcal
  D}E \to M$. Since there is a natural representation $D$ of 
${\mathcal D}F$ on $F$, it follows from Proposition \ref{prop1}
that $\Phi^!(D)$, where $\Phi = \phi^{!!}_{{\mathcal D} F}$, is a
representation of  $\phi^{!!}({\mathcal D} F)$ on $\phi^! F $.
Hence, there is a natural Lie algebroid morphism from
 $ \phi^{!!}  ({\mathcal D}F)$ to ${\mathcal D}(\phi^! F) $.
It suffices to prove that this vector bundle morphism is an isomorphism;
this is an immediate consequence of the fact that the vector bundle 
$\phi^{!!}({\mathcal D} F)$ is locally
isomorphic to the vector bundle
$((F \otimes F^*) \oplus TN) \oplus_{TN} TM$,  which is equal to
$(F \otimes F^*) \oplus TM$, which
in turn is locally isomorphic to $\mathcal D(\phi^! F)$.
\end{proof}
A Lie algebroid morphism $(\Phi, \phi) $ from 
$A \to M$
to $B \to N$ 
is said to be a {\em pull-back morphism}
if $\phi: M \to N $ is admissible and if 
the Lie algebroids $  A \to M   $ and $ \phi^{!!}B \to M $
are isomorphic. In other words, $(\Phi, \phi) $ is a pull-back morphism
when $\Phi$ can be written as 
$$ \Phi = \phi^{!!}_B \circ \Psi \ , $$
where $\Psi  $ is a Lie algebroid isomorphism from $A \to M $ to $
\phi^{!!} B \to M$.
 
\subsection{Pull-backs by transverse maps have vanishing modular classes}

\subsubsection{Transverse maps}
We first recall some well-known facts about short exact sequences.
Let $0 \to V_1 \stackrel{i}{\to} V_2 \stackrel{p}{\to} V_3 \to 0 $ be
a short exact sequence of finite-dimensional vector spaces,
and let $r = {\rm dim} (V_3) $.
There is a canonical isomorphism,
 $$   \wedge^{{\mathrm{top}} } V_1  \otimes \wedge^{{\mathrm{top}}}
 V_3 \stackrel{\simeq}{\to}  
\wedge^{{\mathrm{top}}} V_2 \ ,$$
defined by
$$   X  \otimes  Y  \mapsto  i(X) \wedge \widetilde{Y}  \ ,$$
for all $X \in \wedge^{{\mathrm{top}}}V_1$, $Y \in
\wedge^{{\mathrm{top}}} V_3 $, where $\widetilde{Y}$ is any element in 
$\wedge^{r}  V_2 $ 
such that $(\wedge^{r}p) \, \widetilde{Y} =Y$.
Given an exact sequence $0 \to E_1 \stackrel{i}{\to} E_2
\stackrel{p}{\to} E_3 \to 0 $ 
of vector bundles over a manifold $ M$,  the previous isomorphism 
applied pointwise yields 
a canonical isomorphism of line bundles, 
$$ \wedge^{{\mathrm{top}} } E_1  \otimes \wedge^{{\mathrm{top}}} E_3
\stackrel{\simeq}{\to}  \wedge^{{\mathrm{top}}} E_2 \ ,$$
that we call the {\em canonical line-bundle isomorphism}.

We now define transverse maps.

\begin{definition}\label{transverse}
A map $\phi: M \to N$ is said to be {\em transverse}
to a Lie algebroid $B \to N$ if,  for all $m \in M$, 
 $$ ({T}\phi) ( T_mM  ) \, + \, \rho_B ( B_{\phi(m)} ) = T_{\phi (m) } N .$$
\end{definition}

Requiring 
$\phi $ to be a transverse map amounts to requiring that the following
be an exact sequence of vector bundles,
 \begin{equation}\label{eq:exact_transverse} 0 \to B \oplus_{TN} TM
   \stackrel{i}{\to} 
\phi^{!}B \oplus TM  \stackrel{p}{\to} \phi^{!} TN  \to 0 \ ,
 \end{equation}
where $i$ is the inclusion map, and $p (\phi^{!} b + u ) = \phi^{!}\big(
({T}\phi)\, u -  \rho_B  b \big)$, for all $b \in B$ and $u \in TM$.
In particular, the rank of  
$B \oplus_{TN} TM$ is constant, and the pull-back
Lie algebroid $\phi^{!!} B$ is well-defined. 
We have thus proved the following result.

\begin{lemma}
Any map  $\phi:M \to~N$ transverse to a Lie algebroid $B \to N$  
 is admissible.
\end{lemma}

Let $(\Phi, \phi) $ be a pull-back morphism from 
a Lie algebroid $ A \to M $ to $ B \to N $, with $\phi $
transverse to $B \to N$. By definition, the Lie algebroids $A \to M $ and $  \phi^{!!}B \to M$
are isomorphic. To the short exact sequence (\ref{eq:exact_transverse})
is associated a canonical isomorphism of line bundles,
$$
\wedge^{{\mathrm{top}}} (\phi^{!} TN) \otimes  \wedge^{{\mathrm{top}}} A 
\stackrel{\simeq}{\to}  
  \wedge^{{\mathrm{top}}}  (\phi^{!}B \oplus TM ) = 
 \wedge^{{\mathrm{top}}}  (\phi^{!}B) \otimes  \wedge^{{\mathrm{top}}}
(TM)  \ .
$$
In turn, this isomorphism induces a canonical isomorphism of line
bundles,
\begin{equation} \label{iso}
\wedge^{{\mathrm{top}}} A \otimes
  \wedge^{{\mathrm{top}}}(T^* M)  \stackrel{\simeq}{\to}   
\wedge^{{\mathrm{top}}}
 ( \phi^{!}B)  \otimes \wedge^{{\mathrm{top}}} (\phi^{!} (T^*N) )\ .
\end{equation}
The left-hand side of (\ref{iso}) is the line bundle $L^{A} $ (see
(\ref{linebundle})
in Section \ref{char}), while
the right-hand
side is the line bundle $\phi^{!} L^B $.
In conclusion, in the case of a pull-back 
morphism from $A \to M $ to $B \to N $ over a transverse map, 
(\ref{iso}) is a canonical line bundle isomorphism from $ L^{A} $ to
$\phi^{!} L^B $ that we denote 
by $ \ell^\phi  $. 
In the next sections
we shall prove that 
$$ 
\ell^\phi :  L^{A}   \stackrel{\simeq}{\to}
\phi^{!} L^B 
$$ 
is an isomorphism of $A$-modules.
The proof of the following lemma is straightforward.
\begin{lemma}\label{com_dia}
Let $(\Phi, \phi)$ be a pull-back Lie algebroid morphism from
$A \to M$ to $ B \to N$, with 
$\phi $ transverse to $B \to N$, and let $(\Psi, \psi) $ be a pull-back 
Lie algebroid morphism from $B \to N $ to $C \to R $, with $\psi $
transverse to $C \to R$. Then,

\noindent(i)
$ (\Psi, \psi) \circ (\Phi, \phi)$ is a pull-back Lie algebroid morphism
from $A \to M $  to $C \to R $, with $\psi \circ \phi $ transverse to
$C \to R$, and

\noindent(ii)
the following diagram of line bundle isomorphisms over $M$, 
where $\phi^{!}\ell^{\psi}: \phi^{!}L^B \stackrel{\simeq}{\to}
\phi^{!}\psi^{!} L^C $ is the isomorphism obtained by
pulling back the isomorphism $ \ell^{\psi} : L^B
\stackrel{\simeq}{\to} \psi^{!} L^C $, commutes.
$$
\xymatrix{ L^A  \ar[dr]_{\ell^{\psi \circ \phi}}  \ar[rr]^{\ell^\phi}&
  &  \phi^{!} L^B \ar[dl]^{\phi^{!}\ell^\psi} \\   & \phi^{!}
  \psi^{!}L^C  &   \\}
$$
\end{lemma}
\subsubsection{Pull-backs by submersions}
Any submersion $\phi: M \to N$ is transverse to any Lie
algebroid over $N$.
We shall now prove a property of the pull-backs by submersions.
 \begin{proposition}\label{submersion}
\noindent (i)
When $\phi$ is  a submersion, 
the canonical isomorphism, $\ell^\phi: L^{A} \stackrel{\simeq}{\to}
\phi^{!} L^B $, is
an isomorphism of $A $-modules. 

\noindent(ii)The modular class of a pull-back morphism
$(\Phi, \phi) $ from a Lie algebroid $A \to M$ to a
Lie algebroid $B \to N$ vanishes whenever $\phi: M \to N $
is a submersion.
\end{proposition}
\begin{proof}
By Example \ref{ex:isomorphism} and Equation (\ref{comp}), 
we can assume that $A = \phi^{!!}B $,
and $(\Phi, \phi) = (\phi^{!!}_B , \phi) $, and by 
Theorem \ref{modularmodule}, it suffices to prove
that the canonical isomorphism, 
$\ell^\phi: L^A \stackrel{\simeq}{\to} \phi^{!} L^B $,
is an isomorphism of representations of $ A \to M$.  It suffices, for
that purpose, to prove that $\ell^\phi$
is an isomorphism of representations of $ A_{|_U} \to U$, after
restriction to some neighborhood $U$ of an arbitrary point $m \in M$. 
This restriction allows us to assume that
all the line bundles involved in the computations below 
are topologically trivial, {\it {i.e.}}, admit
global nowhere-vanishing sections.

Let $\mu$ and $\nu$
be volume forms defined on $U \subset M$ and $ \phi(U) \subset 
N$,  respectively, and let $\sigma$ be a nowhere-vanishing section of
$\wedge^{{{\mathrm{top}}}} B$ over $\phi(U)$. 
Henceforth, it will be understood 
that all sections of vector bundles over $M$ (resp., $N$) are 
defined only over $U$ (resp., $\phi(U) $).

Let $F\to M$ be the vertical tangent bundle of the submersion $\phi: M
\to~N$. Then $0 \to F \to TM \to \phi^{!}TN \to 0$ is an exact sequence of
vector bundles over $M$, and there is a canonical isomorphism,
$$ K : \wedge^{{{\mathrm{top}}}} F \otimes \wedge^{{{\mathrm{top}}}}
(T^*M) \stackrel{\simeq}{\to}
\wedge^{{{\mathrm{top}}}} ((\phi^!TN)^*) \stackrel{\simeq}{\to} \phi^!(\wedge^{{{\mathrm{top}}}}
(T^*N)) \ ,
$$  
such that $K(\tau \otimes \mu) = \iota_\tau \mu$, for 
$\tau \in \Gamma(\wedge^{{{\mathrm{top}}}} F)$.
We define $\tau$ to be the unique section of  
$\wedge^{{{\mathrm{top}}}} F$ over $U$
such that $\iota_\tau\mu = \phi^*\nu$.

Since $\phi$ is a submersion, the sequence  $0 \to F \to A \to \phi^!B
\to 0$ is 
exact, so there also exists a canonical isomorphism,
$$
J:  \wedge^{{{\mathrm{top}}}} (\phi^!B) \otimes \wedge^{{{\mathrm{top}}}} F 
\stackrel{\simeq}{\to}
\wedge^{{{\mathrm{top}}}} A \ .
$$
Therefore, one can define a nowhere-vanishing section $ \omega$ of 
$\wedge^{{\mathrm{top}}} A$ over $U \subset M$
by $$
\omega  = J (\phi^! \sigma \otimes \tau ) \ .
$$
It is routine to verify that the canonical isomorphism 
$\ell^{\phi} $ satisfies
\begin{equation}\label{Lphi} \ell^{\phi} \big( \omega  \otimes \mu
  \big)  =    \phi^! (  \sigma \otimes \nu  ) \ . \end{equation}

The modular cocycle for $B$ with respect to the nowhere-vanishing section 
$\sigma \otimes \nu$ of $L^B = \wedge^{\rm{top}} B \otimes
\wedge^{\rm{top}}(T^*N)$ is the section $\beta \in \Gamma(B^*)$ such
that, for any $b \in \Gamma B$,
\begin{equation}\label{I}
<\beta,b> \sigma \otimes \nu = [b,\sigma]_B \otimes \nu + \sigma
\otimes \mathcal L_{\rho_B b}\nu \ .
\end{equation}

The modular cocycle for $\phi^{!!}B$ with respect to the nowhere-vanishing section 
$\omega \otimes \mu$ of
$L^{\phi^{!!}B} = \wedge^{\rm{top}} (\phi^{!!}B) \otimes
\wedge^{\rm{top}}(T^*M)$ is the section $\gamma \in \Gamma((\phi^{!!}B)^*)$ such
that, for any $c \in \Gamma(\phi^{!!} B)$,
\begin{equation}\label{II}
<\gamma,c> \omega \otimes \mu = [c,\omega]_{\phi^{!!}B} \otimes \mu + \omega
\otimes \mathcal L_{u}\mu \ ,
\end{equation}
where $c = \phi^!b + u$, with $b \in \Gamma B$ and $u \in \Gamma(TM)$
such that $\rho_Bb = (T\phi) u$.
 
We shall now prove the relation 
 \begin{equation}\label{eq:goal} \gamma - \widetilde \Phi^*\beta = 0.\end{equation} 
Assuming, without loss of generality, that, locally, $\sigma = b_1 \wedge \ldots \wedge b_p$, with $(b_1,
\ldots, b_p)$ a local basis of sections of $B$, where $p$ is the rank
 of $B$, then $\omega = J(\phi^!\sigma \otimes \tau ) = (\phi^!b_1 + u_1)
\wedge \ldots \wedge  (\phi^!b_p + u_p) \wedge \tau$, where $u_k \in
 \Gamma (TM)$, $k = 1, \ldots, p$. 
Therefore, by  (\ref{bracketforsubmersions}),
$$[c,\omega]_{\phi^{!!}B} = 
[\phi^!b +u, (\phi^!b_1+u_1)\wedge \ldots  \wedge (\phi^!b_p+u_p) \wedge
\tau)]_{\phi^{!!}B}
$$
$$
= \sum_{k=1}^p (-1)^{k+1} [\phi^!b +u, \phi^!b_k+u_k]_{\phi^{!!}B}
\wedge (\phi^!b_1+u_1)\wedge \ldots \wedge (\widehat{\phi^!b_k+u_k})
\wedge \ldots \wedge 
(\phi^!b_p+u_p)\wedge \tau 
$$
$$
+ 
(\phi^!b_1+u_1)\wedge \ldots  \wedge (\phi^!b_p+u_p) \wedge [\phi^!b +u, \tau]_{\phi^{!!}B}
$$
$$
= \! \sum_{k=1}^p (-1)^{k+1} (\phi^![b,b_k]_B  + [u,u_k]_{TM}) 
\wedge (\phi^!b_1+u_1)\wedge \ldots \wedge (\widehat{\phi^!b_k+u_k})
\wedge \ldots \wedge 
(\phi^!b_p+u_p)\wedge \tau 
$$
$$
+
(\phi^!b_1+u_1)\wedge \ldots \wedge  (\phi^!b_p+u_p) \wedge [u, \tau]_{TM} \ .
$$
This formula implies that
$$
[c,\omega]_{\phi^{!!}B} = J\left(\phi^!\sigma \otimes [u,\tau]_{TM} + 
\phi^![b,\sigma]_B \otimes \tau  \right) \ .
$$
Now 
$$
\ell^{\phi}(\omega \otimes \mu) = \ell^{\phi}(J(\phi^!\sigma \otimes \tau) \otimes \mu) =
\phi^!\sigma \otimes \iota_\tau \mu \ = \phi^!\sigma \otimes \phi^* \nu  = 
\phi^*(\sigma \otimes \nu) \ . 
$$
It follows that
$$
<\gamma, c> \ell^{\phi}(\omega \otimes \mu) =
<\gamma,c> \phi^*(\sigma \otimes \nu) 
$$
$$
= \ell^{\phi}\left ( J( \phi^!([b,\sigma]_B \otimes
\tau) \otimes \mu)
+ J (\phi^!\sigma \otimes [u,\tau]_{TM})  \otimes \mu + 
J (\phi^!\sigma \otimes  \tau) \otimes {\mathcal L}_u\mu \right )
$$
$$
= \phi^![b,\sigma]_B \otimes \iota_\tau\mu + \phi^!\sigma \otimes
(\iota_{[u,\tau]_{TM}} \mu + \iota_\tau \mathcal L_u \mu) \ .
$$
From $[\mathcal L_u,\iota_\tau] = \iota_{[u,\tau]}$ and $\iota_\tau\mu
= \phi^*\nu$ 
we obtain
$$
\iota_{[u,\tau]_{TM}} \mu + \iota_\tau \mathcal L_u\mu = \mathcal L_u \iota_\tau \mu =
\mathcal L_u(\phi^* \nu ) = \phi^*(\mathcal L_{(T\phi) u} \nu)
= \phi^*(\mathcal L_{\rho_B b} \nu) \ .
$$
Equation (\ref{eq:goal}) is therefore proved and implies that 
$$ c  \cdot (\omega \otimes \mu)  =\phi^{!} \big(  b \cdot (\sigma
\otimes \nu) \big) \ , $$
for any section $c \in \Gamma A $ of the form $c = \phi^{!} b + u $,
with $\rho_B b =  (T\phi) (u) $. 
This relation, together with  Equations (\ref{Lphi}), (\ref{I}) and
(\ref{II}), 
proves that $\ell^{\phi} $ is an isomorphism of representations of $A \to M $.
\end{proof}

\subsubsection{Pull-backs by transverse maps}
We shall need the following lemma.

\begin{lemma}\label{splitting}
Let $ n \in N$ be a point in the base manifold
of a Lie algebroid $B \to N$,
 $ q$ the rank of the anchor map $\rho_B $ at the point $n  $,
and $ q+s$ the dimension of $N$.
There exist
\begin{itemize}
\item  a surjective submersion $\psi$ from a neighborhood $V$ of $n$ 
to an open subset  $W$ of  $ {\mathbb R}^{s}$, and 
\item a Lie algebroid 
$C \to W $
\end{itemize}
such that
\begin{itemize}
 \item the anchor $\rho_C$ vanishes at the point $\psi (n) $, and
\item the restriction to $V $ of the Lie algebroid $B \to N$ is 
isomorphic to the pull-back $\psi^{!!} C $ of the  Lie algebroid $ C
\to W $ by $\psi$.
\end{itemize}
\end{lemma}
\begin{proof}
By a theorem of Dufour (see Theorem 8.5.1 in
 \cite{DZ}), there exist local coordinates $(x_1, \dots,x_q , y_1, \dots, y_s) $ on $ N$,
centered at $n$, and a local trivialization
 $\alpha_1, \dots, \alpha_q, \beta_1, \dots, \beta_r$ of $B $,
 where $q + r = {\mathrm{rank}} \, B$,  such that, for $1 \leq i, j
 \leq q$ and $1 \leq a, b \leq r$,
 $$
 \left\{ 
\begin{array}{rcl}
 \left[\alpha_i, \alpha_j\right]&=& 0 \\
 \rho_B (\alpha_i) & = & \frac{\partial}{\partial x_i} \\ 
\left[ \beta_a ,\alpha_j \right] &=& 0 \\
 \rho_B(\beta_a) & = & \sum_{k=1}^s g_{a}^k \,  
 \frac{\partial}{\partial y_k} \\
 \left[\beta_a, \beta_b \right] &=& \sum_{c=1}^r
f_{a b}^c \, \beta_c     ,\\
 \end{array}\right. $$
where $g^k_a$ and $f_{ab}^c$ are functions of $y_1, \dots, y_s$, 
and $g^k_a (0, \ldots , 0)=0$.
In particular, the last two 
relations define a Lie algebroid structure $C \to W $ in a neighborhood $W$ of $0 $
in ${\mathbb R}^{s} $ whose anchor vanishes at $0$. 
All together, these formulas show that $B \to N
$ is isomorphic to the pull-back
of $ C \to W$ by the projection $ \psi ( x_1, \dots,x_q , y_1, \dots,
y_s  )  = ( y_1, \dots, y_s) $.
\end{proof}

We can now prove the main result of this section.

\begin{theorem}\label{modtransverse}
The modular class of  a pull-back morphism
$(\Phi, \phi) $ from a Lie algebroid $A \to M$ to a
Lie algebroid $B \to N$, when $\phi: M \to N $ is transverse to $B \to
N$, vanishes.
\end{theorem}
\begin{proof}
As in the proof of Proposition \ref{submersion}, it suffices 
to prove that $\ell^\phi$
is an isomorphism of representations of $ A_{|_U} \to U$, after
restriction to some neighborhood $U$ of an arbitrary point $m \in M $.

By Lemma \ref{splitting}, there exists a neighborhood
$V $ of $n = \phi(m) $, a submersion $\psi: V \to W $,
and a Lie algebroid $ C \to W$ whose anchor vanishes at\break $\psi(n) =
\psi \circ \phi(m)$,
such that $B_{|_V} = \psi^{!!}C $.
Since the anchor map of $C \to W $  vanishes at the point $\psi (n) $,
the differential at the point $m$ of the transverse map 
$\psi \circ \phi $ 
is onto, so that,
shrinking $U$, $V$ and $W$ if necessary, we can assume that $ \psi
\circ \phi$ is a submersion from $U$ onto $W$. 

By Proposition \ref{submersion}, the restriction of 
$\ell^{\psi \circ \phi} $ to $U$ is an isomorphism of representations of 
$A_{|_U} \to U $ from $L^A$ to $(\psi \circ \phi)^!L^C$. Since $\psi$
is a submersion, the restriction of $\ell^{\psi}$ to $V$ also
is an isomorphism of representations of \break
$B_{|_V} \to V $ from $L^B$ to
$\psi^!L^C$. 
Hence, the restriction of $ \phi^{!} \ell^{\psi} $ 
to $U$ is an isomorphism of representations of $A_{|_U} \to U $ from
$\phi^! L^B$ to $\phi^! \psi^! L^C$.
By Lemma \ref{com_dia}(ii), the restriction of $\ell^{\phi} $ to $U $
is also an isomorphism of representations
of $A_{|_U} \to U $ from $L^A$ to $\phi^!L^B$, which completes the proof.
\end{proof}

\begin{proposition}
If $(\Phi, \phi)$ is
a Lie algebroid morphism from $A \to M$ to $B \to N$, and if
$\phi$ is admissible with respect to $B \to N$, the
map $\Phi \colon A \to B$ factors as a composition of maps
$A \stackrel{\Phi'} {\to} \phi^{!!}B \stackrel{\phi^{!!}_B}
{\to} B$, where $\Phi'$ is a morphism of Lie algebroids over $M$.
\end{proposition}

\begin{proof}
When $\phi$ is admissible with respect 
to $B \to N$, the subbundle $\phi^{!!} B$ of
$\phi^!B \oplus TM$ 
is well-defined. For $m \in M$ and $a \in A_m$, the pair $(\phi^!(\Phi(a)),
\rho_A(a))$ is in the fiber of $\phi^{!!} B$ at $m$ 
since $\rho_B \circ \Phi = T\phi
\circ \rho_A$. Setting
$\Phi ' (a) = (\phi^!(\Phi(a)), \rho_A(a))$ defines a  
vector bundle map over the identity of $M$ from $A$ to
$\phi^{!!}B$ such that $\Phi = \phi^{!!}_B \circ \Phi'$.
To prove that $\Phi'$ 
is a Lie algebroid morphism, we write, for any $a \in \Gamma A$, 
$\phi^!(\Phi a) = \sum_i f_i \phi^! (b_i)$, where $f_i \in
C^\infty(M)$
and $b_i \in \Gamma B$, and for $a$ and $a' \in \Gamma A$,
we express the bracket $[ \phi^!(\Phi a), \phi^!(\Phi
a')]_{\phi^{!!}B}$ according to Definition
\ref{pbAlg}. The result then follows from Theorem 4.6 of~\cite{CL}.
\end{proof}

The following is then a corollary of Theorem \ref{modtransverse}.
\begin{proposition} Let $(\Phi, \phi)$ be a Lie algebroid morphism from
  $A\to M $ to $B \to N$ with the map
 $\phi : M \to N$ transverse to $B \to N$, 
and let $\Phi = \phi^{!!}_B \circ \Phi'$ as
 above. Then
$$
{\rm {Mod}} \, {\Phi} = {\rm {Mod}} \, {\Phi'} \ .
$$
\end{proposition}
This result shows that the modular class
of any morphism which covers a transverse map is in fact equal to the
modular class of a base-preserving morphism.

\begin{remark}\label{remarksubmersion}
{\rm {When $\phi$ is a transverse map,  Proposition
\ref{prop1} has the following
simple proof.
The Lie algebroid morphism
$D \colon B \to \mathcal D F$ over $N$ gives rise to a
Lie algebroid morphism $\phi ^{!!} D$ from $\phi^{!!} B$ to $\phi^{!!}
(\mathcal D F)$ over $M$. 
Since, by Proposition \ref{derivations},
$\phi^{!!} (\mathcal D F)$ is isomorphic to $\mathcal D (\phi^! F)$,
the pull-back representation
is obtained as the composition of Lie algebroid morphisms over~$M$, 
$A \stackrel{\Phi'}{\to} \phi^{!!} B \stackrel{\phi^{!!}D}{\to} \mathcal D (\phi^! F)$, a representation of $A$
on $\phi^! F$ satisfying condition \eqref{pullback}.}}
\end{remark}

\subsection{A pull-back morphism with a non-vanishing modular
  class}\label{counterex}
We now give a counter-example which shows that Theorem
\ref{modtransverse} 
does not extend to the case in which $\phi$ is not transverse but only
admissible in the sense of Definition \ref{admissible}.

Let  $N = S^1 \times {\mathbb R} $ be a cylinder
with coordinates $(\theta,x)  $, with $\theta \in S^1 $ and $x \in
{\mathbb R} $.
Let $B \to N$ be the Lie subalgebroid 
of $TN \to N $ generated by the 
vector field   
\begin{equation}\label{eq:defX} {X} = \frac{\partial}{\partial \theta}+ 
x \frac{\partial}{\partial x} ,\end{equation}
whose integral curves spiral in toward the invariant circle $x=0$.
The vector field $X$ is a nowhere-vanishing section
of $\Gamma (\wedge^{\rm{{\mathrm{top}}}}B)=\Gamma B $, while
${\rm d}x \wedge {\rm d}\theta $ is a nowhere-vanishing section
of $\wedge^{\rm{{\mathrm{top}}}} (T^*N) =\wedge^{2} (T^*N) $.
Let $\beta
\in \Gamma (B^* )$ be the modular cocycle for $B$
with respect to the nowhere-vanishing section
${}{X} \otimes ( {\rm d}x \wedge {\rm d}\theta )$
of $ \wedge^{\rm{{\mathrm{top}}}}B
\otimes \wedge^{\rm{{\mathrm{top}}}} (T^*N) $. 
We shall compute 
$\beta$ by evaluating it on $X$,
  $$\begin{array}{lll} <\beta, {}{X} > {}{X} \otimes ({\rm d}x \wedge
    {\rm d}\theta)
& = &  D_{{}{X}}^B  {}\left({X} \otimes ({\rm d}x \wedge {\rm d}\theta)
\right)  \\  
&=&  {}{X} 
\otimes {\mathcal L}_{{}{X}} ({\rm d}x \wedge {\rm d}\theta) 
= {}{X} \otimes  ({\rm d}x \wedge {\rm d}\theta) \ .
\end{array}  $$
In conclusion, the modular class of $B \to N $ is the class of the
section $\beta $ of $B^*$ which satisfies 
\begin{equation}\label{eq:defmodB}<\beta, {}{X} > \, =1 \ .
\end{equation}

Let $M = S^1 $ and let $\phi: M \to N$ be the inclusion map
defined by $\phi(\theta) = (\theta,0)$. We observe that, 
by construction, for all $m \in M $,
\begin{equation}\label{eq:vectorlevel}
   (T \phi) ( \left.
  \frac{\partial}{\partial \theta}\right|_m) = 
  \left. {}{X}\right|_{\phi(m)}.
\end{equation} 
The pair 
$(T \phi , \phi) $ is a Lie algebroid morphism
from $TM \to M $ to $ TN \to N$ 
which takes values in $B \to N$ and yields
 a Lie algebroid morphism  from $ TM \to M$ to $B \to N$,
which we denote also by $(T\phi , \phi) $.
Since the Lie algebroid $TM \to M  $ is unimodular, the modular class
${\rm {Mod}} \, ({ T \phi })$ is $ - \,  (T^*\phi)
({\rm { Mod}} \, B)
$, and the latter is, according to 
 (\ref{eq:defmodB}) and (\ref{eq:vectorlevel}), the class of the
$1$-form  $- \, \beta $ on $S^1$, where $<\beta, \frac{\partial}{\partial
  \theta} > \, =1$. In other words, the relative modular class
of $ (T \phi , \phi)$ is the class of $- \, {\rm d} \theta $, 
and is therefore non-trivial.

Although
the map $\phi $ is not transverse, it is admissible, being the
inclusion of an orbit.  It is easy to show
that the pull-back Lie algebroid
$\phi^{!!} B \to M$ is the Lie algebroid $TM \to M $,
and that the map $ \phi^{!!}: TM \to TN $ is $(T \phi , \phi) $.
Therefore, we have found a Lie algebroid $B$ over $N$ and
an admissible map   $\phi \colon M \to N$
for which the modular class
of the Lie algebroid morphism $\phi^{!!}_B : \phi^{!!} B \to
B $ does not vanish.

\subsection{The $1$-cohomology of the pull-back by a surjective
  submersion}

We first prove a lemma whose result will be used in 
the proof of Proposition \ref{bijective}.

\begin{lemma} \label{globalfunction}
Let $\phi: P \to M $ be a surjective submersion with $1$-connected fibers,
and let $F $ be the vertical tangent bundle of $\phi$.
For any section\break $\beta  \in \Gamma (F^* )$ whose restriction to each fiber
is a closed $1$-form, there exists a smooth function $f \in
C^{\infty}(P)$ such that, for all $u \in F$,
 $${\rm d} f (u) = \beta (u) \ .$$
\end{lemma}
\begin{proof}
Let $(V_i)_{i \in I}$ be an open cover of $M$ which admits local
sections of~$\phi $, 
$ \sigma_i : V_i \to P$, and
let 
$(\chi_i)_{i \in I}$ be a partition of unity associated to $(V_i)$.
On $\phi^{-1}(V_i) $, one can define a smooth function $f_i $ by $f_i (p) = \int_0^1 \beta \big( \gamma'(t) \big) {\rm d}t $,
where $\gamma (t) $ is any path of class $C^1 $ contained in the fiber
of $ p$  with endpoints 
$\gamma (0 ) = \sigma_i (\phi(p))$
and $\gamma (1) =p $.  Since the fibers are $1$-connected, this
function is well-defined. It is differentiable and $ {\rm d} f_i (u)=
\beta (u)$, for any $u \in F_{|V_i} $.
Therefore, the function $ f = \sum_{i \in I} (\chi_i \circ \phi ) \,
f_i $ is a smooth
function on $P$ such that $ {\rm d} f (u) = \beta (u)$, for all $u \in
F$.
\end{proof}
The following proposition (see \cite{C} \cite{F} \cite{GL}) 
describes an important
property of the cohomology of pull-back Lie algebroids.
\begin{proposition} \label{bijective}
Let $A \to M$ be a Lie algebroid, and let
$\phi : P
\to M $ be a surjective submersion with $1$-connected fibers.
The map $\widetilde{(\phi_A^{!!})}^*$ induces an isomorphism from
$H^1(A)$
to $H^1 ( \phi^{!!}A )$.
\end{proposition}
\begin{proof}
We denote $\widetilde{
  (\phi_A^{!!})}^* $  by $\widehat\phi$.
We first prove that the map induced by ${\widehat \phi}$ 
is injective, making no
assumption on the fibers of $\phi$.
Let $\alpha \in \Gamma (A^*) $  be a ${\rm d}_A$-closed section such that
${\widehat \phi} \, \alpha $ is exact, {\it{i.e.}}, such that
${\widehat \phi} \, \alpha = {\rm d}_{\phi^{!!}A} f $, for some function $f \in 
{C}^{\infty} (P)$. We recall that the vertical tangent bundle $F $ 
of the surjective submersion
$ \phi: P \to M$  is a subbundle of $\phi^{!!}A \to P $,
and that $({\widehat \phi}
\alpha) (u) =0$, for all $u  \in F$.
As a consequence, $ {\rm d} f (u) =0 $, for all 
 $u  \in F$, and the function $f$ is basic,
{\it{i.e.}},
there exists a function $g \in  {C}^{\infty} (M) $ such that
$f = g \circ \phi $. Since $(\phi_A^{!!}, \phi)$ is a Lie algebroid morphism,
 $$    {\widehat \phi} \, \alpha={\rm d}_{\phi^{!!}A } f =
{\widehat \phi} ({\rm d}_A g) \  .$$
Since the map $\phi_A^{!!} $ is a surjective bundle map, the above
equation implies that $\alpha = {\rm d}_A g $,
and  therefore that $\alpha  $ is exact. In conclusion, the map induced by
${\widehat \phi}$ in cohomology is injective.

We shall now assume that the fibers of $\phi$ are $1$-connected
in order to prove that the map induced by $\widehat \phi$ is surjective. 
Let $\gamma \in \Gamma ((\phi^{!!} A)^*)  $ be a  
${\rm d}_{\phi^{!!}A}$-closed section. Its restriction to $F $ 
is a fiberwise closed $1$-form that we denote by $\beta$. By Lemma 
\ref{globalfunction}, there exists
a function $f \in {C}^{\infty}(P) $ such that $\beta  (u)= {\rm d}f (u)$ for 
all $u \in F $, {\it{i.e.}}, for all vectors tangent to the fibers of $\phi $.
This implies that  $\gamma ' =  \gamma - {\rm d}_{\phi^{!!}A} f$ is a closed
 section of $ \Gamma ((\phi^{!!} A)^*)$ whose restriction to $ F$ vanishes.
We shall prove that $\gamma' $ is the image under ${\widehat \phi}$ of 
a section $\alpha \in \Gamma (A^* )$.
For $a \in \Gamma A $ and $u \in \Gamma(TP) $ such that $\rho_A a =
(T\phi ) u $, let us set 
$(\phi^{!})^* \alpha (a) = \gamma' (\phi^{!} a + u) $, and let us show
that the $1$-form $\alpha$ on $A$ is well-defined. The 
function $\gamma' (\phi^{!} a + u) $  does not depend on the
choice of $u$ since the restriction of $\gamma'$ to $F $ vanishes. It
is also basic since, spelling out the fact that $\gamma' $ is a
closed section, we obtain,
for any section $v \in \Gamma F$, 
$$ 
0 = ({\rm d}_{\phi^{!!} A}  \gamma ') (\phi^{!}a+u,v) = 
u \cdot \gamma '(v) - v \cdot \gamma '(\phi^{!}a+u)  - \gamma' (
 [\phi^{!}a+u,v]_{\phi^{!!}A} ) \ . $$
Since $v$ and $ [\phi^{!}a+u,v]_{\phi^{!!}A} 
= [u,v]$ are sections of  $F $, 
 $$    v \cdot \gamma '  (\phi^{!}a+u)   =0  \ . $$
By construction, $\alpha $  is a closed section of $A^*$. In
conclusion,  ${\widehat \phi}$ induces a surjective map in cohomology.  
\end{proof}
\section{The case of a regular base-preserving morphism}
\label{fourthsection}

\subsection{The modular class of a totally intransitive Lie algebroid}

Let $A \to M$ be a totally intransitive Lie algebroid, {\it {i.e.}},
such that   
$\rho_A=0$.   Then, for $a \in \Gamma A$ and
$\omega \otimes \mu \in \Gamma(L^A)$, where $L^A = \wedge^{\rm {{\mathrm{top}}}}
A \otimes \wedge^{\rm {top}} (T^*M)$,
$$
D^A_a (\omega \otimes \mu) = 
[a, \omega]_A 
\otimes \mu \ .
$$
By definition, $D^A_a (\omega \otimes \mu) = < \alpha , a > \, \omega
\otimes \mu$, where $\alpha \in \Gamma(A^*)$ is a modular cocycle for
$A$, whence $<\alpha, a > \, \omega = 
[a, \omega]_A$.
Since the anchor of $A$
vanishes, at each point $m \in M$,
$$
< \alpha_{m} , a_m > \, \omega_m = [a_m,\omega_m] \ ,
$$
where the
bracket is the Gerstenhaber bracket on the exterior algebra of the
fiber at $m$ which is a Lie algebra $A_m$. Since $\omega_m$ is a
volume form on $A_m$, 
$\alpha_m \in A_m^*$ is a modular cocycle for the Lie algebra $A_m$,
considered as a Lie algebroid over a point. Therefore
\begin{proposition}\label{prop:LAB}
At each point of a totally intransitive Lie algebroid,
the value of a modular cocycle 
is a modular cocycle for the fiber at that point.
\end{proposition}
In degree $1$, the linear space of cocycles of the Lie algebroid
is the cohomology space $H^1(A)$ of the Lie algebroid, and
similarly for each fiber Lie algebra. Therefore
\begin{corollary}
In a totally intransitive Lie algebroid, the value of the modular class at a point 
is the modular class of the fiber at that point.
\end{corollary}

Whenever $A$ is a Lie algebra bundle,
{\it {i.e.}}, locally trivial as a family of Lie algebras,
 the derived bundle $[A,A]$
of $A$ also is a Lie algebra bundle \cite{M1} \cite{M}, 
and $H^1(A)$ is the space
of sections of the bundle $(A/[A,A])^*$ whose rank is the co-rank of $[A,A]$.
However, when the Lie algebroid $A$ is only totally intransitive,
a Lie algebroid $1$-cocycle on $A$ is a smooth family of Lie
algebra cocycles for the fibers, but
a cocycle for a given fiber might not be extendable to a cocycle on a
neighborhood in $M$, as shown by the example of a
$1$-parameter family of semi-simple algebras 
whose Lie bracket
degenerates to the
zero bracket when the parameter equals zero.  

\subsection{Unimodular extensions of transitive Lie
  algebroids}\label{unimodular}

Let $C$, $A$, $B$ be Lie algebroids over the
base $M$ with anchors $\rho_C$, $\rho_A$ and $\rho_B $,
and brackets $[~,~]_C$, $[~,~]_A$ and $[~,~]_B $,
respectively, and let
$$0 \longrightarrow  
C \stackrel{i}{\longrightarrow}  A
\stackrel{\Phi}{\longrightarrow} 
B \longrightarrow   0  $$
be a Lie algebroid extension \cite{M1} \cite{M}, {\it{i.e.}}, 
$i$ and $\Phi$ are base-preserving Lie algebroid
morphisms, and the sequence is exact. 
Then $C$ is totally intransitive since
$ \rho_C = \rho_A  \circ i = \rho_B \circ \Phi  \circ i
=0$. We assume that $B$ is transitive, which implies that $C$ is a 
Lie algebra bundle (see \cite{M1} \cite{M} for the case $B = TM$).

We shall also assume that the Lie algebra bundle $C$ is unimodular, 
in which case we call the Lie algebroid extension {\it unimodular}.
By Proposition \ref{prop:LAB}, this is the case if and
only if each Lie algebra $C_m$ , $m \in M$, is 
unimodular.  

Finally, for simplicity, we shall assume that the vector
bundles under consideration are orientable, although the results are
valid without this hypothesis.

The adjoint action defines (see \cite{M}) a representation $D^{A,C}$
of the Lie algebroid $A$ on
$C$
by
\begin{equation}
i (D^{A,C}_X k) = [X,i(k)]_A  \ ,
\end{equation}
for $X \in \Gamma A$ and $k \in \Gamma C$, and 
this representation induces a representation $D^{A,K}$ of the Lie
algebroid $A$ on $K = \wedge^{{\rm{top}}} C$ 
such that, for $\lambda \in \Gamma K$, 
\begin{equation}
(\wedge^{r} i ) (D^{A,K}_X \lambda) = 
[X, (\wedge^{r} i) \lambda]_A \ , 
\end{equation}
where $r$ is the rank of $C$ and $ [~,~]_A$
denotes the Gerstenhaber bracket on $\Gamma(\wedge^{\bullet} A)$.
The unimodularity of $C $ implies that $\lambda$ can be chosen so
that, for all $k \in \Gamma C $, $[k,\lambda]_C = 0$. Thus 
$$ 
(\wedge^{r} i ) (D^{A,K}_{i (k)} \lambda) =   
[i (k), (\wedge^{r} i) \lambda]_A  = (\wedge^{r}i)
[k,\lambda]_C =0 \ .
$$
As a consequence, the representation $ D^{A,K} $ factors through the
projection $\Phi : A \to B$, and yields a  representation $D^{B,K}$  of $B$
on $K$ such that 
\begin{equation}
\label{repB}
D^{B,K}_{\Phi (X)} \lambda =  D^{A,K}_X \lambda \ ,
\end{equation}
for all $X \in \Gamma A$.

An immediate consequence of Proposition \ref{prop2}(iii) is

\begin{lemma}\label{lemma1}
The characteristic classes of the representations $ D^{B,K} $
and $D^{A,K}$ are  related by
$$  {\rm{char}} \, {D^{A,K}}  = \Phi^* ({\rm{char}} \, 
{D^{B,K}}) \ .   
$$
\end{lemma}

A computation of the characteristic class of $ D^{A,K}$ yields
 \begin{lemma}\label{lemma2}
The characteristic class of the representation
$D^{A,K}$ is the modular class of $\Phi $.
\end{lemma}
\begin{proof}
Since the sequence $ 0 \to  C_m \stackrel{i}{\longrightarrow}  A_m
\stackrel{\Phi}{\longrightarrow} B_m \to 0$ is exact for all $m \in M $,  
there exists an isomorphism of vector bundles,
$ \check J:  \wedge^{\rm{top}}A 
\otimes  \wedge^{\rm{top}} B^* \to \wedge^{\rm{top}} C   $,
such that
 $$  (\wedge^r i){\check J} (\omega \otimes \mu ) =   \iota_{\Phi^ * \mu} \omega \ ,
$$
for $\omega \in \Gamma(\wedge^{\rm{top}}A)$ and $\mu \in \Gamma
(\wedge^{\rm{top}}B^*)$.

Let $D^\Phi$ be the representation of $A$ on 
$\wedge^{\rm{top}}A \otimes  \wedge^{\rm{top}} B^*$ such that the
characteristic class of $D^\Phi$ is the
modular class ${\rm{Mod}} \, \Phi$ of $\Phi$. Recall from
\cite{KW} that
$$
D^\Phi_X(\omega \otimes \mu) = [X,\omega]_A \otimes \mu + \omega \otimes
{\mathcal L}^B_{\Phi (X)} \mu \ ,
$$
for $X \in \Gamma A$.
The isomorphism $\check J$ intertwines the representations 
$D^\Phi$ and $D^{A,K}$ of $A$. In fact, for any $ X \in \Gamma A $,  
$$( \wedge^ri ) {\check J} (D^\Phi_X ( \omega \otimes \mu ))  =
\iota_{\Phi^* \mu} [X,\omega]_A
+ \iota_{{\mathcal L}^A_X (\Phi^ * \mu)} \omega 
= [X, \iota_{\Phi^* \mu} \omega ]_A $$
$$= [X, (\wedge^r i){\check J}(\omega \otimes \mu)]_A = 
(\wedge^r i) D^{A,K}_X ({\check J} (\omega \otimes \mu)) \ .$$
The conclusion follows since $i$ is injective.
\end{proof}

The preceding lemmas imply the following theorem.

\begin{theorem}\label{th:rel_pb}
For a unimodular Lie algebroid extension, 
$$0 \longrightarrow  
C \stackrel{i}{\longrightarrow}  A
\stackrel{\Phi}{\longrightarrow} 
B \longrightarrow   0  \ ,
$$
the modular class of $\Phi $ satisfies
\begin{equation}\label{extension}
 {\rm{Mod}} \, {\Phi} =\Phi^* ({\rm{char}} \, {D^{B,K}}) \ , 
\end{equation}
where $K = \wedge^{\rm{top}} C$, and $D^{B,K}$ is defined by \eqref{repB}.
\end{theorem}

\begin{remark}{\rm{ It follows from the proofs of Lemmas \ref{lemma1}
      and \ref{lemma2} that the
relation \eqref{extension}
holds at the cochain level.  More precisely, if $\theta$ is the 
$d_A$-cocycle defined by $D^\Phi_X(\omega \otimes \mu) = <\theta,X> \omega \otimes
\mu$, for all $X\in \Gamma A$, and  $\eta$ is the  $d_B$-cocycle
defined by $D^{B,K}_Y \lambda = <\eta, Y>
\lambda$, for all $Y \in \Gamma B$,  then $\theta = \Phi^* \eta$
whenever ${\check J} (\omega
\otimes \mu)= \lambda$.}} 
\end{remark}

\subsection{Morphisms of constant rank with unimodular kernel}\label{unimod}

Let us assume more generally that the  Lie
algebroid morphism
$\Phi \colon A \to A'$ over the identity of $M$ has constant rank,
and that the kernel $C={\rm{ker}} \, \Phi$
is a unimodular Lie algebra bundle. Let $B$ be the image of $A$ under
$\Phi$, and let $\Phi_B$ be
the morphism from $A$ onto $B$ induced by $\Phi$.
To express ${\rm{Mod}} \, {\Phi_B}$, 
we can apply Theorem \ref{th:rel_pb} to the unimodular Lie algebroid extension,
$$
0 \longrightarrow  
C \stackrel{i}{\longrightarrow}  A
\stackrel{\Phi_B}{\longrightarrow} 
B \longrightarrow   0  \ .
$$
Therefore
$${\rm{Mod}} \, A - 
\Phi_B^* ({\rm{Mod}} \, B) =\Phi_B^* 
({\rm{char}} \, {D^{B,K}}) \ .
$$
On the other hand, there is a 
representation $D^{B,\wedge^{\rm{top}}(A'/B)}$ of $B$ on
$\wedge^{\rm{top}}(A'/B)$. If $i_B$ is the 
inclusion of $B$ into $A'$,
then a computation shows that 
$$
i_B^* ({\rm{Mod}} \, A') - {\rm{Mod}} \, B =
{\rm{char}} \, {D^{B,\wedge^{\rm{top}}(A'/B)}} \ .
$$
Since $\Phi = i_B \circ \Phi_B$, and since 
${\rm{Mod}} \, {\Phi} =  {\rm{Mod}} \, A - {\Phi}^*({\rm{Mod}} \, A')$,
we obtain
\begin{theorem}\label{regular}
Let $\Phi : A \to A' $ be a Lie algebroid morphism over the identity
of $M$ whose image $B$ has constant rank and whose kernel 
is a unimodular Lie algebra bundle $C$.  Then the 
modular class of $\Phi $ satisfies
\begin{equation}\label{modregular}
 {\rm{Mod}} \, {\Phi} = \Phi_B^*({\rm{char}} \, {D^{B,K}} -
 {\rm{char}} \, {D^{B,\wedge^{\rm{top}}(A'/B)}}) \ ,
\end{equation}
where $K = \wedge^{\rm{top}} C$, and $D^{B,K}$ is defined by \eqref{repB}.
\end{theorem}

\subsection{Examples}

\subsubsection{Transitive Lie algebroids}\label{transitive}
We apply the results of Section \ref{unimod} to the case
where $B=TM$ and 
$\Phi$ is the anchor map $\rho_A$ of a transitive Lie algebroid $A$. The
 modular class of $\Phi$ is then the
modular class of the Lie algebroid $A$. 
As in Section \ref{unimodular}, we assume that the isotropy bundle $C =
{\rm{ker}} \rho_A$ is unimodular.

In this case, the representation $D^{B,K}$ of 
$B=TM$  
is a natural flat connection $\nabla$ on $K = \wedge^{\rm{top}}C$. 
The characteristic class
${\rm{char}} \, {D^{B,K}}$
is the class of the closed $1$-form $\alpha$ on $M$ such that, for all
$X \in \Gamma(TM)$ and $\lambda \in \Gamma K$, 
$\nabla_X \lambda = <\alpha, X > \lambda$.
Therefore Theorem
\ref{th:rel_pb} implies that
$${\rm{Mod}} \, A = \rho_A^* ({\rm{char}}(D^{TM,\wedge^{\rm{top}}(\rm{ker} \rho_A)}) \ .$$
In particular, if $M$ is simply-connected, any transitive Lie
algebroid over $M$ 
with unimodular isotropy bundle is unimodular.

\subsubsection{Regular Poisson structures} Let $E$ be a Lie
algebroid with a Poisson or twisted Poisson structure \cite{KL}~\cite{KY} 
defined by a bivector $\pi \in \Gamma(\wedge^2E)$. Assume that the
structure is regular, 
{\it {i.e.}}, 
the associated map, $\pi^{\sharp} : E^* \to E$, is of
constant rank. Applying Theorem \ref{th:rel_pb}
to $A=E^*$, $B=
{\rm {Im}}(\pi^{\sharp})$ and $\Phi = \pi^\sharp_B$, the submersion
from $E^*$ to $B= {\rm {Im}}(\pi^{\sharp})$ defined by $\pi^\sharp$, 
we recover Lemma 2.3 of \cite{KY}, which is written, in the notations
that we are using here,
$$
{\rm{Mod}} \, {\pi^\sharp_B} = 
(\pi^\sharp_B)^*({\rm{char}} \, {D^{B,K}}) \ . $$
In this case, the representation 
$D^{B,\wedge^{\rm{top}}(A'/B)}$ is dual to $D^{B,K}$, and therefore
Theorem \ref{regular} yields
\begin{equation}\label{poisson}{\rm{Mod}} \, {\pi^\sharp} = 2 \,
(\pi^\sharp_B)^* ({\rm{char}} \, {D^{B,K}}) \ , 
\end{equation}
which is Theorem 2.5 of \cite{KY}. 

In the case of a Poisson manifold,  $(M,\pi)$, the integrable
distribution $B \subset TM$ defines the symplectic foliation, and 
\eqref{poisson} implies  Corollary 9 of \cite{C}. In particular (see \cite{W}),
\begin{proposition}
A regular Poisson manifold is unimodular if and only if its 
symplectic foliation 
admits an invariant transverse volume form.
\end{proposition}

\section{Generalized morphisms of Lie algebroids}\label{fifthsection}

\subsection{The category of Lie algebroids with 
generalized Lie algebroid morphisms}
There are several ways to define what could be 
called generalized morphisms of Lie algebroids,
as an analogue to the definition of the generalized morphisms of Lie 
groupoids \cite{BuW} \cite{MM}.  Here, we follow the basic idea of
Ginzburg \cite{Gi}, who defined Morita equivalences, but allowing
enough generality to include all ordinary morphisms.

\begin{definition}
A generalized morphism ${\mathcal P} $ from a Lie algebroid $A \to
M$ to a Lie algebroid $B \to N$
is a manifold $P$ together with 

\noindent{$\bullet$}
a surjective submersion $\phi$ from $P $ to
$M$, and

\noindent{$\bullet$} a Lie algebroid morphism $(\Psi,\psi) $ from $\phi^{!!}A \to
P$ to $B \to N $. 
\end{definition}
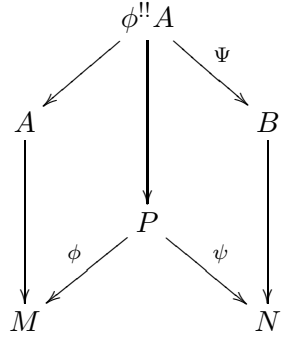
\vspace{-.4cm}
$$
\xymatrix{  & \ar[dl]_{} \phi^{!!} A \ar[dd] \ar[dr]^{\Psi}  \\ A
  \ar[dd]{} &  &  B \ar[dd] \\
    & \ar[dl]_{\phi} P \ar[dr]^{\psi} & \\
M & & N    }
$$

We denote a generalized morphism by a triple, ${\mathcal P} = 
(P ,\phi, (\Psi,\psi))$.

An {\em isomorphism} between generalized morphisms ${\mathcal P}_i
=(P_i, \phi_i, (\Psi_i. \psi_i))$, $i=1,2 $, from
$A \to M$ to $B \to N$
is  a diffeomorphism $\sigma: P_1 \buildrel{\simeq}\over{\to} P_2 $
that intertwines both $\phi_1 $ and $\phi_2 $, and $\Psi_1 $ and
$\Psi_2 $. More precisely, it is  a diffeomorphism such 
that the following diagrams commute.
$$ 
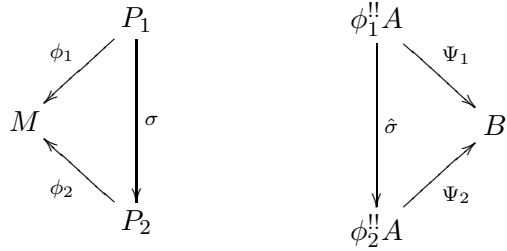
\xymatrix{  & \ar[dl]_{\phi_1} P_1 \ar[dd]^{\sigma}  \\ M  &    \\
  &  \ar[ul]^{\phi_2} P_2   } \quad \quad \quad  \quad \quad \quad 
\xymatrix{   \phi_1^{!!} A  \ar[dr]^{\Psi_1}  \ar[dd]^{\hat \sigma}
  & \\ & B  \\  \phi_2^{!!} A \ar[ur]_{\Psi_2}  &  }  $$
where ${\hat \sigma}$ 
is the Lie algebroid isomorphism from  $\phi_1^{!!} A$ to
$\phi_2^{!!} A$
induced by $\sigma $, defined by 
$${\hat \sigma} ( \phi_1^{!} a+ u ) = \phi_2^{!}a \, + \,  (T\sigma )
\, u \ ,$$
 for all $a \in A_{\phi_1 (p)}$, $p \in P_1$, $u \in T_p P_1 $, 
where $ \rho_A a = (T \phi_1) u $. It is clear that 
the isomorphism of generalized morphisms is an equivalence relation.

We shall now prove that generalized morphisms can be composed.
Let  ${\mathcal P} = (P, \phi, (\Psi, \psi)) $ be a generalized morphism from 
a Lie algebroid $A \to M$ to a Lie algebroid $B \to N$, and let
${\mathcal P}' = (P', \phi', (\Psi', \psi')) $ be a generalized morphism from 
a Lie algebroid $B \to N$ to a Lie algebroid $C \to R$.
Consider the following data:

\noindent{$\bullet$}
the set $P'' = P \times_{\psi , N, \phi'} P' $, which is a manifold
since $\phi' $ is a surjective submersion from $P'$ to $N$,
 
\noindent{$\bullet$}
the map $\phi''$ defined by 
$\phi ''( p, p'  ) = \phi (p) $ for all $(p,p') \, \in \,  P'' $, 
which is a surjective submersion,

\noindent{$\bullet$}
the Lie algebroid $ (\phi'')^{!!} A \to P''$ whose 
fiber over $ (p,p') \in P''$ is 
$$
\left\{ (\phi'')^{!} a +u + v  \in    
(\phi'')^{!} A_p \oplus T_p P \oplus T_{p'} P' \, \mid \,   \rho_A a =
(T\phi) u \, , \, (T\psi)u = (T\phi')v  \right\}. 
$$
We can then define a vector bundle map $\Psi''$ from $(\phi'')^{!!} A$
to $C$ by
  $$\Psi'' ((\phi'')^{!} a + u + v ) = 
\Psi' \left( (\phi')^{!}\Psi (\phi^{!}a + u) + v \right) \  .$$
This vector bundle map is, in fact, a Lie algebroid
morphism covering the map $\psi'' : P \to R$ defined by
$\psi''(p,p') = \psi'(p')$. 
Therefore ${\mathcal P}''=(P'', \phi'', (\Psi'',\psi''))$ is
a generalized morphism 
from $A \to M$ to  $C \to R$.
The compositions of isomorphic generalized morphisms are isomorphic,
and the associativity of the composition is valid up to isomorphism.
These properties justify
the following definition.
\begin{definition} (i)
The category of Lie algebroids with generalized Lie algebroid morphisms 
is the
category $\algd '$ in which
\begin{itemize}
\item objects are Lie algebroids, 
\item arrows are isomorphism classes of generalized morphisms,
and the identity arrow 
of a Lie  algebroid $A \to M$
is the isomorphism class of the 
generalized morphism $(M, {\rm id}_M, ({\rm id}_A,{\rm id}_M))$.
\end{itemize}

\noindent(ii)
A {\em Morita equivalence} from a Lie algebroid $A \to M$ to a Lie algebroid
$B \to N$ is a generalized morphism
${\mathcal P}=(P,\phi, (\Psi, \psi)) $ 
where $( \Psi, \psi) $ is a pull-back Lie algebroid morphism covering 
a surjective submersion $\psi : P \to N$.
\end{definition}

\subsection{Modular classes of generalized morphisms}
We now consider generalized Lie algebroid morphisms
${\mathcal P}=(P, \phi, (\Psi, \psi )) $ such that the fibers of the
map $\phi: P \to M$ are $1$-connected. 
The composition of two such generalized morphisms satisfies
the same condition (see \cite{Gi}), and therefore  
generalized Lie algebroid morphisms 
with 
$1$-connected fibers
define a sub-category $\algd'_1$ of the category $\algd'$
of Lie algebroids with generalized Lie
algebroid morphisms.

We recall from
Proposition \ref{bijective} that a pull-back morphism covering a
surjective submersion $\phi$ with $1$-connected fibers induces an 
isomorphism of Lie algebroid cohomologies in degree~$1$, 
which we have denoted by
$\widehat \phi$.
We can now define the modular class of a generalized morphism with
$1$-connected fibers.

\begin{definition}
When ${\mathcal P}= (P,\phi, (\Psi, \psi)) $ 
is a Lie algebroid generalized morphism with $1$-connected fibers, its 
modular class is the image of the modular class of $\Psi$
by ${\widehat \phi}^{-1}$.
\end{definition}
The modular class of 
$ {\mathcal P}$ is an element of $H^1(A)$ which we denote by
${\rm Mod } \, {\mathcal P}$.
When ${\mathcal P}= (P,\phi, (\Psi, \psi)) $ is a generalized morphism 
from $A \to M$ to $B \to N$
such that the fibers of $\phi$ are $1$-connected, we denote by
${\widetilde{\mathcal P}}^* $
the map
from $H^1 (B) $ to $H^1 (A) $
induced by $\widehat{\phi}^{-1} \circ \widetilde \Psi^*$.
It is clear from diagram (\ref{eq:comdia_iso}) below that isomorphic
generalized morphisms induce the 
same map on the $1$-cohomology.

\begin{theorem}\label{composition}

\noindent(i) The modular classes of  
isomorphic generalized morphisms with $1$-connected fibers are equal.

\noindent(ii) The modular class of the composition of two generalized
morphisms ${\mathcal P} $ and ${\mathcal P}' $
with $1$-connected fibers is 
\begin{equation}\label{composed}
 {\rm Mod} ({\mathcal P}' \circ {\mathcal P}) =
{\rm Mod} \, {\mathcal P} + \widetilde{{\mathcal P}}^* ( {\rm Mod}
 \, {\mathcal P}') \ .
\end{equation}

\noindent(iii) The modular class of a Morita equivalence with $1$-connected
   fibers vanishes.
\end{theorem}

\begin{proof}
(i) Let $\sigma: P_1 \to P_2$ be an isomorphism between generalized 
morphisms ${\mathcal P}_i = (P_i, \phi_i,(\Psi_i, \psi_i)) $, $i=1,2$, from
a Lie algebroid  $A \to M$ to a Lie algebroid $B \to N$.
Denote by $\hat{\sigma} $  the Lie algebroid isomorphism from
$\phi_1^{!!}A$ to $\phi_2^{!!}A  $  induced by $\sigma $.
The following diagram commutes. 
\begin{equation}
\label{eq:comdia_iso} \xymatrix{   &  \phi_1^{!!} A \ar[dr]^{\Psi_1}
  \ar[dd]^{\hat{\sigma}} 
\ar[dl]_{(\phi_1)_A^{!!}} &   \\ A  &   & B \\   & \phi_2^{!!}A \ar[ur]_{\Psi_2} 
\ar[ul]^{(\phi_2)_A^{!!}} &   } \end{equation}
Since $\hat \sigma $ induces an isomorphism of
Lie algebroid cohomologies in degree~$1$, and since, by Proposition
\ref{bijective},
the modular classes of $\phi_1^{!!} A$ and $\phi_2^{!!} A$
vanish, the result follows from the commutativity of the above diagram.

(ii)
The projection onto the
first component, $\chi \! : \! P'' \! = \! P \times_{\psi , N ,
  \phi'}P' \to~P $, 
is  a surjective submersion, and there is a natural identification
of Lie algebroids $\chi^{!!} (\phi^{!!}A ) \stackrel{\simeq}{\to}
 (\phi'')^{!!} A$.
There is also 
a Lie algebroid morphism $\Xi: (\phi'')^{!!} A \to (\phi')^{!!} B  $, 
covering the projection onto the second component $\xi :P'' = 
P \times_{\psi , N . \phi'} P' \to P' $,
defined by $$ \Xi \big( (\phi'')^{!} a +u + v \big)= (\phi')^{!} \Psi (
\phi^{!} a +u) + v \  ,$$ for $a \in A$, $u \in TP$, $v \in TP'$
such that  $\rho_A a = (T\phi) u$ and
$(T\psi)u = (T\phi')v$, and 
the following  diagram of Lie algebroid morphisms commutes.
$$
\xymatrix{ &  & (\phi'')^{!!} A  \ar@/^2pc/[rrdd]^{\Psi''} \ar@/_2pc/[lldd]_{( \phi'')^{!!}_A } 
 \ar[dl]^{\chi^{!!}_{\phi^{!!}A}}  \ar[dr]_{\Xi}& & \\ &\phi^{!!} A  \ar[dl]^{\phi^{!!}_A}  \ar[dr]_{\Psi}& &(\phi')^{!!} B  \ar[dl]^{(\phi')^{!!}_B}  \ar[dr]_{\Psi'} & \\
           A &  & B &  & C \\ }
$$
By Proposition \ref{bijective}, all arrows pointing ``south-west''
induce isomorphisms of Lie algebroid cohomologies  
in degree~$1$. Since, by Proposition
\ref{submersion}, each has a vanishing modular class,
the result follows from the commutativity of the above diagram.

(iii) is an immediate consequence of Proposition \ref{submersion}.
\end{proof}

Theorem \ref{composition} (i) shows that the modular class of an isomorphism
class of generalized morphisms with $1$-connected fibers is well
defined,
while (\ref{composed}) generalizes formula (\ref{comp}).

\section{Appendix: Representations of categories and the 
modular class} \label{Appendix}

\subsection{Representations of categories}
Let $\bf R$ be a commutative ring.
A {\em representation} of a category $\mathcal C$ over $\bf R$ is a functor $F$ from
$\mathcal C$ to the
category of $\bf R$-modules \cite{BW} \cite{G}.
An {\em anti-representation} of a category $\mathcal C$ is a
representation of the opposite category,
 $\mathcal C^{\rm {opp}}$.

When the category $\mathcal C$ is a Lie groupoid 
$G\arrows X$, and $\bf R$ is the field of real numbers, we assume that the
vector spaces $F(x)$ associated to the objects, {\it {i.e.}}, points, $x$ in
$X$ are the fibers of a smooth vector bundle.  The assignment of vector
bundle morphisms to the morphisms in $G$ is just an action of $G$ in
the usual sense; we assume that this action is smooth.

\begin{example}\rm{
A representation of a Lie group considered as a groupoid
over a point is a representation in the usual sense.}
\end{example}

\begin{example}\rm{
A representation of the action groupoid
associated to the action of a group $H$ on a manifold $X$ is the same
as an $H$-equivariant vector bundle over $X$.}
\end{example}

\begin{remark}\rm{
Representations of Lie algebroids are ``infinitesimal
representations''; they may occur as linearizations of representations
of Lie groupoids.}
\end{remark}

\subsection{$H^1$ as an anti-representation of the 
category of Lie algebroids}

Assigning to each object $A$ of  $\mathcal{A}{{\rm{lgd}}}$ 
the vector space $H^1(A)$ and to each Lie algebroid morphism $\Phi : A
  \to B$, where $(A,B)$ is a pair of objects in ${\mathcal
  A}{\rm{lgd}}$, 
the linear map
  $\widetilde{\Phi}^* : H^1(B) \to H^1(A)$
defines an anti-representation of this category over $\mathbb R$, 
{\it{i.e.}}, a functor from ${\mathcal A}{\rm{lgd}}^{\rm{opp}}$ to the category of 
$\mathbb R$-modules, which we may denote by $H^1$.  The results of
Section \ref{fifthsection} show that $H^1$ is also an anti-representation of
$\algd'_1$.  

\subsection{Cohomology of a category with values in an
  anti-represen\-tation.}

The $0$-cochains on a category $\mathcal C$
with values in an anti-representation $F$ 
assign to any object $A$ of $\mathcal C$ 
an element of $F(A)$. 
The $k$-cochains, $k \geq 1$, 
are maps from $k$-tuples of composable morphisms with source $A$, where
$A$
is an object of $\mathcal C$,
to $F(A)$.
The differential $\delta$ is defined by the usual formula.
This definition is simpler than the one which we have found
in the literature
(see, {\it{e.g.}}, \cite{BW} \cite{G}), 
but it defines the same object in the cases studied here,
and seems therefore to be good enough for our limited purposes.

In particular, the coboundary of a $0$-cochain $u$ with values in $F$ 
is the $1$-cochain $\delta u$
defined by
$(\delta u)(\Phi) = u(A) - F(\Phi)(u (B))$, for any
morphism\break $ \Phi : A \to B$.
If $v$ is a $1$-cochain with values in $F$, its coboundary $\delta v$ 
assigns to 
a pair $(\Phi, \Psi)$ of composable morphisms,
$A \stackrel{\Phi}{\rightarrow} B \stackrel{\Psi}{\rightarrow} C$, where $(A,B,C)$ is
a triple
of objects in the category, the element  
$(\delta v) (\Phi, \Psi)$ which is the alternating sum
$v (\Phi) - v(\Psi \circ \Phi) +
F(\Phi)(v(\Psi))$ in
$F(A)$.

\subsection{ The modular class of morphisms 
as a $1$-coboundary.}

A $0$-cochain on $\mathcal{A}{\rm{lgd}}$ with values in $H^1$
assigns to any object $A$ of $\mathcal{A}{\rm{lgd}}$ an element of $H^1(A)$.
The modular class of Lie algebroids is a $0$-cochain ${\rm{Mod}}$
on $\mathcal{A}{\rm{lgd}}$ (resp.,  on ${\algd'_1}$) with values in 
the anti-representation $H^1$.  

A $1$-cochain  assigns to a morphism $\Phi : A \to B$ an element of $H^1(A)$.
The coboundary of a $0$-cochain $u$ on $\mathcal{A}{\rm{lgd}}$ with values in 
$H^1$ is the $1$-cochain $\delta u$
defined by
$(\delta u)(\Phi) = u(A) - {\widetilde{\Phi}}^* (u (B))$. In particular,
$$
(\delta {\rm{Mod}})(\Phi) = {\rm{Mod}} \, A - {\widetilde{\Phi}}^*
({\rm{Mod}} \, B) \ .$$   The same argument holds for generalized morphisms. 
Thus, we have the following result.
\begin{proposition} 
The modular class of morphisms (resp., generalized morphisms with $1$-connected fibers) of Lie algebroids
is the coboundary of the modular class of Lie algebroids, viewed as a
  $0$-cochain 
on $\mathcal{A}{\rm{lgd}}$ (resp., on
${\algd'_1}$) with values in the anti-representation $H^1$.

\end{proposition} 

If $v$ is a $1$-cochain, its coboundary $\delta v$ assigns to 
a composable pair\break 
$A \stackrel{\Phi}{\to} B \stackrel{\Psi}{\to} C$ of Lie algebroid
morphisms the element in $H^1(A)$, 
$$
(\delta v) (\Phi, \Psi) = v (\Phi) - v(\Psi \circ \Phi) +
{\widetilde{\Phi}}^*(v(\Psi)) \ .
$$
As a consequence of the relation $\delta^2= 0$, 
if $v = \delta u$, 
\begin{equation}\label{eq3}
v(\Psi \circ \Phi) = v (\Phi) 
+ {\widetilde{\Phi}}^*(v(\Psi)) \ .
\end{equation}
The preceding relation
applied 
to $v = \delta {\rm{Mod}}$ yields 
Equations (\ref{comp}) and (\ref{composed}).

\medskip

\begin{remark}
{\rm
One may ask whether there are any characteristic classes of Lie
algebroid morphisms which satisfy \eqref{eq3} but which do
not arise from characteristic classes of the objects in the category
of Lie algebroids.  
In fact, there
are none, for the following general reason.
In any category $\mathcal C$ with an initial object $\{{\rm{pt}}\}$,
{\it {i.e.}}, 
the opposite of a category with a final object, such as 
the category of Lie algebroids, the cohomology in degree $1$ is
trivial. In fact, for each object $A$ of $\mathcal C$, let us denote
by $p_A$ the morphism $A \rightarrow \{{\rm{pt}}\}$, and let $v$ be a
$1$-cochain on $\mathcal C$. Let $u$ be the $0$-cochain defined by
$u(A) = v(p_A)$.  If $v$ is a $1$-cocycle, for any morphism $\Phi : A
\to B$, $v(p_B) = v(p_A \circ \Phi) = v(\Phi) + F(\Phi) (v(p_A))$,
whence $v(\Phi) = u(B) - F(\Phi) (u(A))$, and $v$ is the coboundary of
$u$.

On the other hand, there do exist nontrivial cocycles on other
interesting categories.   We learned the following example from Peter
Teichner.  Let $\mathcal C$ be the category whose objects are smooth
manifolds and whose morphisms are immersions $f:M\to N$ equipped with
an almost complex structure on the normal bundle.  (Composition of
morphisms is composition of immersions and direct sum of complex
structures.)  The de Rham $1$-cohomology is a representation of this
category, and the Chern class of the normal bundle is a $1$-cocycle
which is not a coboundary.
}
\end{remark}

\bigskip

\medskip

\noindent
      {\it Yvette Kosmann-Schwarzbach}\\
      Centre de Math\'ematiques Laurent Schwartz \\
      \'Ecole Polytechnique\\
      91128 Palaiseau, France\\
      yks@math.polytechnique.fr

\medskip

\noindent
{\it Camille Laurent-Gengoux} \\
D{\'e}partement de Math{\'e}matiques\\
Universit{\'e} de Poitiers \\
86962 Futuroscope Cedex, France\\
Camille.Laurent@math.univ-poitiers.fr

\medskip

\noindent
{\it  Alan Weinstein }\\
 Department of Mathematics \\
 University of California \\
        Berkeley, CA 94720-3840, USA \\
 alanw@math.berkeley.edu  \\
\end{document}